\documentclass[10pt, oneside, a4paper]{article}
\usepackage{amsmath, amsthm, amssymb, amsfonts}
\usepackage{geometry}
\usepackage{graphicx}
\usepackage{natbib}
\usepackage{setspace}
\usepackage{algorithm}
\usepackage{longtable}
\usepackage{caption}
\usepackage{xr}

\newtheorem{theorem}{Theorem}[section]
\newtheorem{lemma}[theorem]{Lemma}

\theoremstyle{definition}

\newtheorem{assumption}{Assumption}

\input{tcilatex}

\begin{document}

\title{A Necessary and Sufficient Condition for Size Controllability of
Heteroskedasticity Robust Test Statistics\thanks{%
We thank Mikkel S\o lvsten for helpful discussions and for suggesting to
re-express condition (\ref{non-incl_Het_uncorr}) as condition (\ref{sol}) in
Remark 2.1(ii). We are also grateful to two referees and a Co-Editor for
helpful comments.}}
\author{ 
\begin{tabular}{c}
Benedikt M. P\"{o}tscher\thanks{%
Corresponding author.} \\ 
{\footnotesize University of Vienna} \\ 
{\footnotesize Department of Statistics} \\ 
{\footnotesize A-1090 Vienna, Oskar-Morgenstern Platz 1} \\ 
{\footnotesize benedikt.poetscher@univie.ac.at}%
\end{tabular}
\and 
\begin{tabular}{c}
David Preinerstorfer \\ 
{\footnotesize WU Vienna University of Economics and Business} \\ 
{\footnotesize Institute for Statistics and Mathematics} \\ 
{\footnotesize A-1020 Vienna, Welthandelsplatz 1} \\ 
{\footnotesize david.preinerstorfer@wu.ac.at}%
\end{tabular}
}
\date{First version: June 2024\\
Second version: August 2024\\
Third version: November 2024\\
Fourth version: September 2025\\
Fifth version: April 2026\\
This version: June 2026}
\maketitle

\begin{abstract}
We revisit size controllability results in \cite{PP21} concerning
heteroskedasticity robust test statistics in regression models. For the
special, but important, case of testing a single restriction (e.g., a zero
restriction on a single coefficient), we povide a necessary \emph{and}
sufficient condition for size controllability, whereas the condition in \cite%
{PP21} is, in general, only sufficient (even in the case of testing a single
restriction).
\end{abstract}

\section{Introduction\label{Intro}}

Tests and confidence intervals based on so-called heteroskedasticity robust
standard errors date back to \cite{E63, E67} and constitute, at least since 
\cite{W80}, a major component of the applied econometrician's toolbox.
Although these early methods come with well-understood large sample
properties, when based on critical values derived from asymptotic theory
their finite sample properties often deviate substantially from what
asymptotic theory suggests: tests may substantially overreject under the
null and corresponding confidence intervals may undercover. Strong leverage
points have been identified early on as one major reason for these
deviations, see, e.g., \cite{MacW85}, \cite{DavidsonMacKinnon1985}, and \cite%
{CheshJewitt1987}. This has led to various developments trying to attenuate
such drawbacks:

\begin{enumerate}
\item modifications of the covariance matrix estimators in \cite{E63, E67}
and \cite{W80} led to tests based on what are now frequently called HC1-HC4
covariance matrix estimators (see, e.g., \cite{LE2000}, and \cite{Crib2004}
for an overview of the relevant literature), with HC0 denoting the original
proposal;

\item some authors investigated degree-of-freedom corrections to obtain
modified critical values (e.g., \cite{Satterth} or \cite{BellMcCa}, see also 
\cite{Imbkoles2016});

\item wild bootstrap methods were investigated (for an overview of the
relevant literature see \cite{PPBoot}) and, more recently, parametric
bootstrap methods were studied in \cite{Chuetal2021} and \cite{Hansen2021}.
\end{enumerate}

Although these developments sometimes lead to improvements, they come with
no general finite sample guarantees concerning the size of the tests or the
coverage of related confidence intervals, cf.~the discussion in \cite%
{PPBoot,PP21} for detailed accounts.

Motivated by this lack of finite sample guarantees, \cite{PP21} studied the
question under which conditions heteroskedasticity robust test statistics as
well as the standard (uncorrected) F-test statistic can actually be paired
with appropriate (finite) critical values, so that one obtains tests that
have their (finite sample) size controlled by the prescribed significance
value $\alpha $ (i.e., have size $\leq \alpha $) even though one is
completely agnostic about the form of heteroskedasticity.\footnote{%
The null-hypothesis to be tested is given by a set of affine restrictions.}
Under appropriate assumptions on the errors, allowing for Gaussian as well
as substantial non-Gaussian behavior, they have shown that the standard
(uncorrected) F-test statistic can be size-controlled (in finite samples) by
using an appropriately chosen (finite) critical value if and only if the
following simple condition holds: 
\begin{eqnarray}
&&\text{\emph{no standard basis vector that lies in the column span of the
design matrix}}  \notag \\
&&\text{\emph{is \textquotedblleft involved\textquotedblright\ in the affine
restrictions to be tested,}}  \label{q:cond}
\end{eqnarray}%
see (8) in \cite{PP21} for a formal statement of this condition.

Under a generally \emph{stronger} condition than (\ref{q:cond}) (see (10) in 
\cite{PP21}), it was furthermore shown that large classes of
heteroskedasticity robust test statistics (e.g., HC0-HC4) can be
size-controlled by appropriate (finite) critical values. That condition,
however, although satisfied for many testing problems (and even often
identical to (\ref{q:cond}), cf.~Theorem 3.9 and Lemma A.3 in \cite{PP3}),
is \emph{not} necessary in general, as shown in examples given in \cite{PP21}%
; e.g., their Example 5.5 or Example C.1 in their Appendix C.\footnote{%
Appendices to \cite{PP21} are published in the Supplementary Material
available at the publisher's website of that article.} These examples
consider the case of testing linear contrasts in the expected outcomes of
subjects belonging to two or more groups, scenarios that are practically
relevant. Further examples are provided in Examples A.1-A.4 in Appendix \ref%
{app:proofs} further below.\footnote{%
Example 5.5 in \cite{PP21} concerns simultaneously testing multiple
retrictions, while Example C.1 in Appendix C of \cite{PP21} as well as
Examples A.1-A.4 in Appendix \ref{app:proofs} of the present article concern
the case of testing a single restriction.}

For the important case of testing problems involving only a \emph{single
restriction }(i.e., the case $q=1$ in the notation of \cite{PP21}), we show
in the present article that the condition in (\ref{q:cond}) is then in fact
necessary and sufficient also for size controllability of the above
mentioned classes of heteroskedasticity robust test statistics, including
HC0-HC4.

\section{Results on size controllability\label{frame}}

\subsection{Framework\label{frm}}

Here we recall the most relevant notions from Sections 2 and 3 of \cite{PP21}%
, to which we refer the reader for further information and discussion. We
consider the linear regression model 
\begin{equation}
\mathbf{Y}=X\beta +\mathbf{U},  \label{lm}
\end{equation}%
where $X$ is a (real) nonstochastic regressor (design) matrix of dimension $%
n\times k$ and where $\beta \in \mathbb{R}^{k}$ denotes the unknown
regression parameter vector. Throughout, we assume $\limfunc{rank}(X)=k$ and 
$1\leq k<n$. We furthermore assume that the $n\times 1$ disturbance vector $%
\mathbf{U}=(\mathbf{u}_{1},\ldots ,\mathbf{u}_{n})^{\prime }$ ($^{\prime }$
denoting transposition) has mean zero and unknown covariance matrix $\sigma
^{2}\Sigma $ ($0<\sigma <\infty $), where $\Sigma $ varies in the
\textquotedblleft heteroskedasticity model\textquotedblright\ given by%
\begin{equation*}
\mathfrak{C}_{Het}=\left\{ \limfunc{diag}(\tau _{1}^{2},\ldots ,\tau
_{n}^{2}):\tau _{i}^{2}>0\text{ for all }i\text{, }\sum_{i=1}^{n}\tau
_{i}^{2}=1\right\} ,
\end{equation*}%
and where $\limfunc{diag}(\tau _{1}^{2},\ldots ,\tau _{n}^{2})$ denotes the
diagonal $n\times n$ matrix with diagonal elements given by $\tau _{i}^{2}$.
That is, the disturbances are uncorrelated but can be heteroskedastic of
arbitrary form. [In Appendix \ref{app:B} we shall also consider another
heteroskedasticity model.]\footnote{%
Since we are concerned with finite-sample results only, the elements of $%
\mathbf{Y}$, $X$, and $\mathbf{U}$ (and even the probability space
supporting $\mathbf{Y}$ and $\mathbf{U}$) may depend on sample size $n$, but
this will not be expressed in the notation. Furthermore, the obvious
dependence of $\mathfrak{C}_{Het}$ on $n$ will also not be shown in the
notation, and the same applies to the heteroskedasticity model defined in
Appendix \ref{app:B}.}

\emph{For ease of exposition, we shall maintain in the sequel that the
disturbance vector }$\mathbf{U}$ \emph{is normally distributed.
Generalizations to classes of non-normal disturbances can be obtained
following the arguments in Section 7.1 of \cite{PP21}, see Remark 2.2
further below.} Denoting a Gaussian probability measure with mean $\mu \in 
\mathbb{R}^{n}$ and (possibly singular) covariance matrix $A$ by $P_{\mu ,A}$%
, the collection of distributions on $\mathbb{R}^{n}$ (the sample space of $%
\mathbf{Y}$) induced by the linear model just described together with the
Gaussianity assumption is then given by 
\begin{equation*}
\left\{ P_{\mu ,\sigma ^{2}\Sigma }:\mu \in \mathrm{\limfunc{span}}%
(X),0<\sigma ^{2}<\infty ,\Sigma \in \mathfrak{C}_{Het}\right\} ,
\end{equation*}%
where $\mathrm{\limfunc{span}}(X)$ denotes the column space of $X$.\footnote{%
Since every $\Sigma \in \mathfrak{C}_{Het}$ is positive definite, the
measure $P_{\mu ,\sigma ^{2}\Sigma }$ is absolutely continuous with respect
to Lebesgue measure on $\mathbb{R}^{n}$.}

We focus on testing the null $R\beta =r$ against the alternative $R\beta
\neq r$, where $R\neq 0$ is a $1\times k$ vector and $r\in \mathbb{R}$. That
is, \emph{throughout this paper we focus on testing a single restriction,
whereas the theory developed in \cite{PP21} allows for simultaneously
testing multiple restrictions} (that is, we here consider only the special
case corresponding to $q=1$ in \cite{PP21}). Set $\mathfrak{M}=\limfunc{span}%
(X)$, define the affine space 
\begin{equation*}
\mathfrak{M}_{0}=\left\{ \mu \in \mathfrak{M}:\mu =X\beta \text{ and }R\beta
=r\right\} ,
\end{equation*}%
and let 
\begin{equation*}
\mathfrak{M}_{1}=\left\{ \mu \in \mathfrak{M}:\mu =X\beta \text{ and }R\beta
\neq r\right\} .
\end{equation*}%
Adopting these definitions, the testing problem we consider can be written
more precisely as 
\begin{equation}
H_{0}:\mu \in \mathfrak{M}_{0},\ 0<\sigma ^{2}<\infty ,\ \Sigma \in 
\mathfrak{C}_{Het}\quad \text{ vs. }\quad H_{1}:\mu \in \mathfrak{M}_{1},\
0<\sigma ^{2}<\infty ,\ \Sigma \in \mathfrak{C}_{Het}.
\label{testing problem}
\end{equation}%
We also write~$\mathfrak{M}_{0}^{lin}=\mathfrak{M}_{0}-\mu _{0}=\left\{
X\beta :R\beta =0\right\} $ where $\mu _{0}\in \mathfrak{M}_{0}$. Of course, 
$\mathfrak{M}_{0}^{lin}$ does not depend on the choice of $\mu _{0}\in 
\mathfrak{M}_{0}$. Furthermore, if $\mathcal{L}$ is a linear subspace of $%
\mathbb{R}^{n}$, $\Pi _{\mathcal{L}}$ denotes the orthogonal projection onto 
$\mathcal{L}$, while $\mathcal{L}^{\bot }$ denotes the orthogonal complement
of $\mathcal{L}$ in $\mathbb{R}^{n}$.

The assumption of nonstochastic regressors made above entails little loss of
generality, and results for models with stochastic regressors can be
obtained from the ones derived in the present paper by the same arguments as
the ones given in Section 7.2 of \cite{PP21}.

\subsection{Test statistics, size controllability, and a new result\label%
{Sec_2.2}}

We consider the same test statistics as in Section 3 of \cite{PP21}.
Simplified to the setting of testing a \emph{single} restriction considered
in the present article, they are given by 
\begin{equation}
T_{Het}\left( y\right) =\left\{ 
\begin{array}{cc}
(R\hat{\beta}\left( y\right) -r)^{2}/\hat{\Omega}_{Het}(y) & \text{if }\hat{%
\Omega}_{Het}\left( y\right) \neq 0, \\ 
0 & \text{if }\hat{\Omega}_{Het}\left( y\right) =0,%
\end{array}%
\right.  \label{T_het}
\end{equation}%
where $\hat{\beta}(y)=\left( X^{\prime }X\right) ^{-1}X^{\prime }y$ and
where $\hat{\Omega}_{Het}(y)=R\hat{\Psi}_{Het}(y)R^{\prime }$. Here%
\begin{equation*}
\hat{\Psi}_{Het}\left( y\right) =(X^{\prime }X)^{-1}X^{\prime }\limfunc{diag}%
\left( d_{1}\hat{u}_{1}^{2}\left( y\right) ,\ldots ,d_{n}\hat{u}%
_{n}^{2}\left( y\right) \right) X(X^{\prime }X)^{-1},
\end{equation*}%
with $\hat{u}(y)=\left( \hat{u}_{1}(y),\ldots ,\hat{u}_{n}(y)\right)
^{\prime }=y-X\hat{\beta}(y)$. The constants $d_{i}>0$ sometimes depend on
the design matrix; see \cite{PP21} for examples of the weights $d_{i}$,
including HC0-HC4 weights. We also recall the following assumption from the
latter reference, again specialized to the setting of testing only a \emph{%
single} restriction (i.e., to the case $q=1$ in the notation of \cite{PP21}).

\begin{assumption}
\label{R_and_X}Let $1\leq i_{1}<\ldots <i_{s}\leq n$ denote all the indices
for which $e_{i_{j}}(n)\in \limfunc{span}(X)$ holds where $e_{j}(n)$ denotes
the $j$-th standard basis vector in $\mathbb{R}^{n}$. If no such index
exists, set $s=0$. Let $X^{\prime }\left( \lnot (i_{1},\ldots i_{s})\right) $
denote the matrix which is obtained from $X^{\prime }$ by deleting all
columns with indices $i_{j}$, $1\leq i_{1}<\ldots <i_{s}\leq n$ (if $s=0$,
no column is deleted). Then $R(X^{\prime }X)^{-1}X^{\prime }\left( \lnot
(i_{1},\ldots i_{s})\right) \neq 0$ holds.
\end{assumption}

This assumption can be checked in any particular application as it only
depends on the observable quantities $R$ and $X$; and a sufficient condition
for Assumption \ref{R_and_X} obviously is $s=0$.\footnote{%
To check Assumption \ref{R_and_X}, one first determines the indices $i_{j}$ (%
$j=1,\ldots ,s$), which amounts to checking which of the matrices $%
(X:e_{i}(n))$ have rank equal to $k$. Given these indices, one then computes 
$R(X^{\prime }X)^{-1}X^{\prime }\left( \lnot (i_{1},\ldots i_{s})\right) $,
and checks whether or not this is zero. We note that checking Assumption \ref%
{R_and_X} is implemented via the sub-routine \textquotedblleft
As.check\textquotedblright\ in the R package \textbf{hrt}, which accompanies 
\cite{PP21}, cf. also the discussion in Appendix E there. The subroutine
\textquotedblleft As.check\textquotedblright\ can be found in the package's
\textquotedblleft auxiliary.functions.R\textquotedblright\ file, to which we
refer the interested reader for details (e.g., the precise algorithm used to
numerically determine the rank of a matrix).} Assumption \ref{R_and_X} is
unavoidable if one wants to obtain a sensible test from the statistic $%
T_{Het}$, see Section 3 of \cite{PP21} for more discussion. We note that $%
e_{j}(n)\in \limfunc{span}(X)$ is equivalent to $h_{jj}=1$, where $h_{jj}$
denotes the $j$-th diagonal element of the `hat matrix' $H=X(X^{\prime
}X)^{-1}X^{\prime }$.\footnote{%
This follows from $h_{jj}=e_{j}(n)^{\prime }He_{j}(n)=(He_{j}(n))^{\prime
}He_{j}(n)$ and the fact that $H$ represents the orthogonal projection onto $%
\limfunc{span}(X)$.}

As in \cite{PP21}, we introduce%
\begin{equation*}
B(y)=R(X^{\prime }X)^{-1}X^{\prime }\limfunc{diag}\left( \hat{u}%
_{1}(y),\ldots ,\hat{u}_{n}(y)\right) .
\end{equation*}%
Define (recall that $R$ is a nonzero row vector in this article) 
\begin{equation*}
\mathsf{B}=\left\{ y\in \mathbb{R}^{n}:\limfunc{rank}(B(y))<1\right\}
=\left\{ y\in \mathbb{R}^{n}:B(y)=0\right\} .
\end{equation*}%
It is now easy to see that $\limfunc{span}(X)\subseteq \mathsf{B}$ and that $%
\mathsf{B}$ is a linear space (cf.~also Lemma 3.1 in \cite{PP21}). Simple
examples can be constructed to show that $\limfunc{span}(X)\neq \mathsf{B}$,
in general; cf.~Example C.1 in Appendix C of \cite{PP21} as well as Examples
A.1-A.4 in Appendix \ref{app:proofs} further below.

To summarize the main size controllability statements from \cite{PP21} for
the above class of test statistics, we first have to recall the following
notation: For a given linear subspace $\mathcal{L}$ of $\mathbb{R}^{n}$ we
define the set of indices $I_{0}(\mathcal{L})$ via 
\begin{equation}
I_{0}(\mathcal{L})=\left\{ i:1\leq i\leq n,e_{i}(n)\in \mathcal{L}\right\} .
\label{eqn:I0def}
\end{equation}%
We set $I_{1}(\mathcal{L})=\left\{ 1,\ldots ,n\right\} \backslash I_{0}(%
\mathcal{L})$. Clearly, $\func{card}(I_{0}(\mathcal{L}))\leq \dim (\mathcal{L%
})$ holds. And $I_{1}(\mathcal{L})$ is nonempty provided $\dim (\mathcal{L}%
)<n$; in particular, $I_{1}(\mathfrak{M}_{0}^{lin})$ is always nonempty
since $\dim (\mathfrak{M}_{0}^{lin})=k-1<n-1$. The results in \cite{PP21}
concerning size controllability of tests for (\ref{testing problem}) based
on $T_{Het}$ can now be summarized as follows; some intuition for why size
control cannot always be achieved is provided further below as well as in
Section 4 in \cite{PP21}:

\begin{theorem}[Theorem 5.1(b,c) and Propositions 5.5(b) and 5.7(b) in 
\protect\cite{PP21} for the case $q=1$]
\label{thm:pp21}\footnote{%
The corresponding results in \cite{PP21} for $q\geq 1$ take exactly the same
form, but with the definitions of the relevant quantities adapted to that
more general setting.} Suppose that Assumption \ref{R_and_X} is satisfied.
Then the following statements hold:

\begin{enumerate}
\item For every $0<\alpha <1$ there exists a real number $C(\alpha )$ such
that 
\begin{equation}
\sup_{\mu _{0}\in \mathfrak{M}_{0}}\sup_{0<\sigma ^{2}<\infty }\sup_{\Sigma
\in \mathfrak{C}_{Het}}P_{\mu _{0},\sigma ^{2}\Sigma }(T_{Het}\geq C(\alpha
))\leq \alpha  \label{size-control_Het}
\end{equation}%
holds, provided that 
\begin{equation}
e_{i}(n)\notin \mathsf{B}\text{ \ \ for every \ }i\in I_{1}(\mathfrak{M}%
_{0}^{lin}).  \label{non-incl_Het}
\end{equation}%
Furthermore, under condition (\ref{non-incl_Het}), even equality can be
achieved in (\ref{size-control_Het}) by a proper choice of $C(\alpha )$,
provided $\alpha \in (0,\alpha ^{\ast }]\cap (0,1)$ holds, where 
\begin{equation}
\alpha ^{\ast }=\sup_{C\in (C^{\ast },\infty )}\sup_{\Sigma \in \mathfrak{C}%
_{Het}}P_{\mu _{0},\Sigma }(T_{Het}\geq C)  \label{alpha}
\end{equation}%
is positive and where 
\begin{equation}
C^{\ast }=\max \{T_{Het}(\mu _{0}+e_{i}(n)):i\in I_{1}(\mathfrak{M}%
_{0}^{lin})\}  \label{Cstar}
\end{equation}%
for $\mu _{0}\in \mathfrak{M}_{0}$ (with neither $\alpha ^{\ast }$ nor $%
C^{\ast }$ depending on the choice of $\mu _{0}\in \mathfrak{M}_{0}$).

\item Suppose (\ref{non-incl_Het}) is satisfied. Then a smallest critical
value, denoted by $C_{\Diamond }(\alpha )$, satisfying (\ref%
{size-control_Het}) exists for every $0<\alpha <1$. And $C_{\Diamond
}(\alpha )$ is also the smallest among the critical values leading to
equality in (\ref{size-control_Het}) whenever such critical values exist.

\item Suppose (\ref{non-incl_Het}) is satisfied. Then any $C(\alpha )$
satisfying (\ref{size-control_Het}) necessarily has to satisfy $C(\alpha
)\geq C^{\ast }$. In fact, for any $C<C^{\ast }$ we have $\sup_{\Sigma \in 
\mathfrak{C}_{Het}}P_{\mu _{0},\sigma ^{2}\Sigma }(T_{Het}\geq C)=1$ for
every $\mu _{0}\in \mathfrak{M}_{0}$ and every $\sigma ^{2}\in (0,\infty )$.

\item If the condition%
\begin{equation}
e_{i}(n)\notin \func{span}(X)\text{ \ \ for every \ }i\in I_{1}(\mathfrak{M}%
_{0}^{lin})  \label{non-incl_Het_uncorr}
\end{equation}%
is violated, then $\sup_{\Sigma \in \mathfrak{C}_{Het}}P_{\mu _{0},\sigma
^{2}\Sigma }(T_{Het}\geq C)=1$ for \emph{every} choice of critical value $C$%
, every $\mu _{0}\in \mathfrak{M}_{0}$, and every $\sigma ^{2}\in (0,\infty
) $ (implying that size equals $1$ for every $C$).\footnote{\label{FN}It is
understood here that critical values are less than infinity.}
\end{enumerate}
\end{theorem}

To obtain some intuition for Theorem \ref{thm:pp21}, recall that the
diagonal elements of $\Sigma \in \mathfrak{C}_{Het}$ are positive and sum up
to one (by definition). Now, for a matrix $\Sigma $ with $i$-th diagonal
entry close to $1$, all other diagonal entries must therefore be close to $0$%
, so that $\Sigma \approx e_{i}(n)e_{i}(n)^{\prime }$ then holds. Note that
if $\Sigma \approx e_{i}(n)e_{i}(n)^{\prime }$, the distribution $P_{\mu
,\sigma ^{2}\Sigma }$ of the data is strongly \textquotedblleft
concentrated\textquotedblright\ around the one-dimensional space $\mu +%
\limfunc{span}(e_{i}(n))$. From an intuitive point of view, whether a given
test statistic admits a size-controlling critical value or not, should
therefore depend on the \textquotedblleft behavior\textquotedblright\ of the
test statistic for values on or close to the spaces $\mu _{0}+\limfunc{span}%
(e_{i}(n))$ with $\mu _{0}\in \mathfrak{M}_{0}$. It turns out that this is
intimately related to (\ref{non-incl_Het}) and (\ref{non-incl_Het_uncorr}).
See Section 4 in \cite{PP21} for more discussion.

Most importantly, the above theorem shows that, given Assumption~\ref%
{R_and_X}, the condition in (\ref{non-incl_Het}) is sufficient for the
existence of a (finite) size-controlling critical value $C(\alpha )$
satisfying (\ref{size-control_Het}), while the weaker condition (\ref%
{non-incl_Het_uncorr}) is necessary. Furthermore, in case the design matrix $%
X$ and the vector $R$ are such that $\mathsf{B}=\limfunc{span}(X)$, and
hence the condition in (\ref{non-incl_Het}) coincides with that in (\ref%
{non-incl_Het_uncorr}), the condition (\ref{non-incl_Het}) is also
necessary. However, $\mathsf{B}=\limfunc{span}(X)$ is not always true (see
Example C.1 in Appendix C of \cite{PP21} or the examples in Appendix \ref%
{app:proofs} further below), although the equality holds generically
(cf.~Theorem 3.9 and Lemma A.3 in \cite{PP3}). We now show in the subsequent
theorem that in the situation considered in this article, namely testing
only a single restriction, the condition in (\ref{non-incl_Het}) in Theorem %
\ref{thm:pp21} can actually always be replaced by that in (\ref%
{non-incl_Het_uncorr}). Before we present that theorem, we discuss an
equivalent formulation of condition (\ref{non-incl_Het_uncorr}) that is
expressed in terms of \emph{certain} diagonal elements of the `hat matrix' $%
H $, see (\ref{sol}) below.\footnote{%
An informal verbal description of (\ref{non-incl_Het_uncorr}) is given in (%
\ref{q:cond})\ in the Introduction.}

\textbf{Remark 2.1:} (i) Condition (\ref{non-incl_Het_uncorr}) is equivalent
to "$h_{ii}<1$ for every $i\in I_{1}(\mathfrak{M}_{0}^{lin})$".\footnote{%
Note that $h_{ii}=1$ always holds if $i\in I_{0}(\mathfrak{M}_{0}^{lin})$.}

(ii) Condition (\ref{non-incl_Het_uncorr}) can also equivalently be written
as 
\begin{equation}
e_{i}(n)\notin \func{span}(X)\text{ for every }i\text{ satisfying }%
R(X^{\prime }X)^{-1}x_{i\cdot }^{\prime }\neq 0,  \label{sol0}
\end{equation}%
see Remark B.1(iii) in Appendix \ref{app:B} further below.\footnote{%
Comparing (\ref{non-incl_Het_uncorr}) and (\ref{sol0}) could lead one to
conjecture equivalence of the conditions $i\in I_{1}(\mathfrak{M}_{0}^{lin})$
and $R(X^{\prime }X)^{-1}x_{i\cdot }^{\prime }\neq 0$. This is incorrect in
general, see Example A.1 in Appendix A. However, $R(X^{\prime
}X)^{-1}x_{i\cdot }^{\prime }\neq 0$ implies $i\in I_{1}(\mathfrak{M}%
_{0}^{lin})$, see Part 3 of Lemma \ref{gensol}.} And this in turn is now
equivalent to%
\begin{equation}
h_{ii}<1\text{ for every }i\text{ satisfying }R(X^{\prime }X)^{-1}x_{i\cdot
}^{\prime }\neq 0.  \label{sol}
\end{equation}%
The last form of the condition may be more appealing to some readers. We
issue a warning here, however, namely that the condition (\ref{non-incl_Het}%
) is, in general, stronger than the condition "$e_{i}(n)\notin \mathsf{B}$
for every $i$ satisfying $R(X^{\prime }X)^{-1}x_{i\cdot }^{\prime }\neq 0$",
see Remark B.1(iv) in Appendix \ref{app:B}.

We now present the announced theorem.

\begin{theorem}
\label{thm:pp24}\footnote{%
Following the suggestion of some readers, we mention the equivalent
condition (\ref{sol}) explicitly in this theorem, although the equivalence
with (\ref{non-incl_Het_uncorr}) has already been noted in Remark 2.1.}
Suppose that Assumption \ref{R_and_X} is satisfied. Then the following
statements hold:

\begin{enumerate}
\item For every $0<\alpha <1$ there exists a real number $C(\alpha )$ such
that (\ref{size-control_Het}) holds, provided that (\ref{non-incl_Het_uncorr}%
) (or equivalently (\ref{sol})) holds. Furthermore, under condition (\ref%
{non-incl_Het_uncorr}) (or equivalently (\ref{sol})), even equality can be
achieved in (\ref{size-control_Het}) by a proper choice of $C(\alpha )$,
provided $\alpha \in (0,\alpha ^{\ast }]\cap (0,1)$ holds, where $\alpha
^{\ast }$ given by (\ref{alpha}) is positive and where $C^{\ast }$is given
by (\ref{Cstar}) for $\mu _{0}\in \mathfrak{M}_{0}$ (with neither $\alpha
^{\ast }$ nor $C^{\ast }$ depending on the choice of $\mu _{0}\in \mathfrak{M%
}_{0}$).

\item Suppose (\ref{non-incl_Het_uncorr}) (or equivalently (\ref{sol})) is
satisfied. Then a smallest critical value, denoted by $C_{\Diamond }(\alpha
) $, satisfying (\ref{size-control_Het}) exists for every $0<\alpha <1$. And 
$C_{\Diamond }(\alpha )$ is also the smallest among the critical values
leading to equality in (\ref{size-control_Het}) whenever such critical
values exist.

\item Suppose (\ref{non-incl_Het_uncorr}) (or equivalently (\ref{sol})) is
satisfied. Then any $C(\alpha )$ satisfying (\ref{size-control_Het})
necessarily has to satisfy $C(\alpha )\geq C^{\ast }$. In fact, for any $%
C<C^{\ast }$ we have $\sup_{\Sigma \in \mathfrak{C}_{Het}}P_{\mu _{0},\sigma
^{2}\Sigma }(T_{Het}\geq C)=1$ for every $\mu _{0}\in \mathfrak{M}_{0}$ and
every $\sigma ^{2}\in (0,\infty )$.

\item If (\ref{non-incl_Het_uncorr}) (or equivalently (\ref{sol})) is
violated, then $\sup_{\Sigma \in \mathfrak{C}_{Het}}P_{\mu _{0},\sigma
^{2}\Sigma }(T_{Het}\geq C)=1$ for \emph{every} choice of critical value $C$%
, every $\mu _{0}\in \mathfrak{M}_{0}$, and every $\sigma ^{2}\in (0,\infty
) $ (implying that size equals $1$ for every $C$).\footnote{%
Cf. Footnote \ref{FN}.}
\end{enumerate}
\end{theorem}

The main take-away of Theorem \ref{thm:pp24} is that, given Assumption \ref%
{R_and_X} holds, the condition in (\ref{non-incl_Het_uncorr}) (or
equivalently (\ref{sol})) is necessary and sufficient for the existence of a
(smallest) finite size-controlling critical value when one is testing only a
single restriction.\footnote{%
By contraposition, given Assumption \ref{R_and_X} holds, the design matrices 
$X$ and restrictions $R$ for which size control fails are precisely
characterized by failure of\ (\ref{non-incl_Het_uncorr}) (or equivalently (%
\ref{sol})). One example of failure of\ (\ref{non-incl_Het_uncorr}) (or
equivalently (\ref{sol})) is when $X$ contains the dummy $e_{i}(n)$ as its
first column, say, and $R=(1,0,\ldots ,0)$ (or, more generally, $R$ has a
non-zero first entry). Another example arises when the first two columns of $%
X$ are given by $(1,\ldots ,1)^{\prime }$ and $(1,-1,\ldots ,-1)^{\prime }$,
and the first two entries of $R$ are both equal to $1$. [Instances of these
two examples for which Assumption \ref{R_and_X} is also satisfied are easily
found.]} The condition "$e_{i}(n)\notin \func{span}(X)$ for every $%
i=1,\ldots ,n$" (which is tantamount to "$h_{ii}<1$ for every $i=1,\ldots ,n$%
") implies (\ref{non-incl_Het_uncorr}), and thus is sufficient for
size-controllability of $T_{Het}$ (but not necessary, see, e.g., Example
A.2). Note that the conditions in (\ref{non-incl_Het}), (\ref%
{non-incl_Het_uncorr}), as well as (\ref{sol}) do not depend on the weights
used in the construction of the covariance matrix estimator or on $r$. They
only depend on $X$ and $R$. This and more (e.g., how the conditions relate
to high-leverage points) is discussed subsequent to Theorem 5.1 (and in
Remarks 5.2-5.4, 5.6, and 5.9) in \cite{PP21} to which we refer the reader
for a detailed account. As a point of interest we also note that condition (%
\ref{non-incl_Het_uncorr}) given above is exactly the same as condition (8)
in \cite{PP21} (with $q=1$); in that reference, the latter condition is
shown to be necessary and sufficient for size control of the standard
(uncorrected) F-test statistic (regardless of whether $q=1$ or not).

We also note here that Theorem \ref{thm:pp24} disproves -- for the special
case of testing a single restriction -- a conjecture in Remark 5.8 of \cite%
{PP21}, namely that there would exist cases where Assumption \ref{R_and_X}
holds, (\ref{non-incl_Het_uncorr}) is satisfied, (\ref{non-incl_Het}) does
not hold, and size control by a (finite) critical value is not possible.

To see why the refinement of Theorem \ref{thm:pp21} provided in Theorem \ref%
{thm:pp24} can matter in practice, it is enough to consider the textbook
example of a matrix $X$ with two columns, the first indicating membership to
the treatment group and the second indicating membership to the control
group (a special case of Example C.1 in Appendix C of \cite{PP21}). Assume
that the first $n_{1}\geq 2$ observations belong to the treatment group and
the remaining $n_{2}\geq 2$ observations belong to the control group.%
\footnote{%
In case $n_{1}=1$ and $n_{2}\geq 2$, condition (\ref{non-incl_Het_uncorr})
is violated. However, this is of no interest, as in this case Assumption \ref%
{R_and_X} is also violated, making Theorem \ref{thm:pp24} inapplicable. In
fact, in this case the test statistic $T_{Het}$ is identically zero, and
hence is size controllable in a trivial fashion. Cf. the discussion in
Example C.1 in Appendix C of \cite{PP21}.} Assume further that one wants to
test whether $\beta _{1}$, the expected outcome of the treatment group,
equals a given value (e.g., because one wants to obtain a confidence
interval through test inversion). Example C.1 in Appendix C of \cite{PP21}
shows that in this case 
\begin{equation*}
\mathcal{I}_{1}(\mathfrak{M}_{0}^{lin})=\{1,\ldots ,n\}\quad \text{ and }%
\quad \mathsf{B}=\{y\in \mathbb{R}^{n}:y_{1}=\ldots =y_{n_{1}}\}\neq \func{%
span}(X).
\end{equation*}%
In particular, $e_{i}(n)\in \mathsf{B}$ if and only if $i>n_{1}$, so that (%
\ref{non-incl_Het}) is not satisfied, while (\ref{non-incl_Het_uncorr})
holds, and size-controlling critical values hence exist by Theorem \ref%
{thm:pp24} (and can be used for constructing confidence intervals). Further
examples are provided in Appendix \ref{app:proofs} below.

We next explain the key observation underlying the proof of Theorem \ref%
{thm:pp24}: To this end, define the (possibly empty) set of indices 
\begin{equation*}
\mathcal{I}_{\#}=\left\{ i:1\leq i\leq n,~R(X^{\prime }X)^{-1}x_{i\cdot
}^{\prime }=0\right\} ,
\end{equation*}%
where $x_{i\cdot }$ denotes the $i$-th row of $X$, and define (the span of
the empty set will throughout be interpreted as $\{0\}$) the space 
\begin{equation}
\mathcal{V}_{\#}=\limfunc{span}\left( \{e_{i}(n):i\in \mathcal{I}%
_{\#},~e_{i}(n)\in \mathsf{B}\}\right) \subseteq \mathsf{B},
\label{eqn:Vdef}
\end{equation}%
the inclusion holding because $\mathsf{B}$ is a linear space as noted
earlier (recall that $R$ is $1\times k$ dimensional in this article).%
\footnote{%
We note that $\mathcal{I}_{\#}$ is a proper subset of $\{1,\ldots ,n\}$
since $R\neq 0$.} Recall that under Assumption \ref{R_and_X} the test
statistic $T_{Het}$ as well as $\mathsf{B}$ are invariant with respect to
(w.r.t.) the group $G(\mathfrak{M}_{0})$ (i.e., the group of transformations 
$y\mapsto \delta (y-\mu _{0})+\mu _{0}^{\ast }$ with $\delta \in \mathbb{R}$
nonzero and $\mu _{0}$ and $\mu _{0}^{\ast }$ in $\mathfrak{M}_{0}$), see
Remark C.1 in Appendix C of \cite{PP21}.\footnote{%
The invariance holds trivially if Assumption \ref{R_and_X} is violated.} The
results in \cite{PP21} are based on this invariance property. The crucial
observation exploited in the proof of Theorem \ref{thm:pp24} now is that, in
the special case of testing a single restriction considered in this article,
the test statistic $T_{Het}$ as well as $\mathsf{B}$ are invariant, not only
w.r.t.~$G(\mathfrak{M}_{0})$, but also w.r.t.~addition of elements of $%
\mathcal{V}_{\#}$. This additional invariance property involving $\mathcal{V}%
_{\#}$, paired with a careful application of the general theory for
size-controlling critical values in \cite{PP3}, then allows us to deduce the
refined statement in Theorem \ref{thm:pp24}. It turns out fortunate that the
general theory in \cite{PP3} explicitly allows one to incorporate additional
invariance properties beyond $G(\mathfrak{M}_{0})$. For details and proofs
the reader is referred to Appendices \ref{app:proofs} and \ref{app:B}.

Finally, we remark that Theorem \ref{thm:pp24} is deduced from Theorem \ref%
{thm:moregen} in Appendix \ref{app:B}, which is a more general statement
that also allows for heteroskedasticity models other than $\mathfrak{C}%
_{Het} $ (and which are defined in (\ref{eqn:covmodmg}) below).

\bigskip

\textbf{Remark 2.2:} \emph{(Extensions to non-Gaussian errors) }(i) All the
theorems in this article continue to hold as they stand, if the disturbance
vector $\mathbf{U}$ follows an elliptically symmetric distribution that has
no atom at the origin; more precisely, $\mathbf{U}$ is assumed to be
distributed as $\sigma \Sigma ^{1/2}\mathbf{z}$, where $\mathbf{z}$ has a
spherically symmetric distribution on $\mathbb{R}^{n}$ that has no atom at
the origin, and where $\sigma $ and $\Sigma $ are as in Section \ref{frm}.
This is so, since the size under Gaussianity is the same as the size under
the elliptical symmetry assumption. In particular, the smallest
size-controlling critical values under the elliptical symmetry assumption
coincide with the smallest size-controlling critical values under
Gaussianity, and thus can be computed from the algorithms relying on
Gaussianity described in \cite{PP21}. See Appendix E.1 of \cite{PP3} and
Section 7.1(i) of \cite{PP21} for more details. The same is actually true
for a wider class of distribution for $\mathbf{U}$, namely where $\mathbf{z}$
has a distribution in the class $Z_{ua}$ defined in Appendix E.1 of \cite%
{PP3}.

(ii) All the theorems in this article except for Theorem \ref{thm:moregen}
in Appendix \ref{app:B} (i.e., all theorems using the heteroskedasticity
model $\mathfrak{C}_{Het}$) continue to hold as they stand, if it is assumed
that the disturbance vector $\mathbf{U}$ follows a distribution from the
semiparametric model defined in Section 7.1(iv) in \cite{PP21} (a model that
contains inter alia all distributions corresponding to i.i.d. samples of
scale-mixtures of normals). Again, this is so since the size under
Gaussianity is the same as the size under this semiparametric model. In
particular, the smallest size-controlling critical values under this
semiparametric model coincide with the smallest size-controlling critical
values under Gaussianity, and thus can be computed from the algorithms
relying on Gaussianity described in \cite{PP21}. See Section 7.1(iv) in \cite%
{PP21} and note that the Gaussian model is a submodel of the semiparametric
model considered there.

(iii) Furthermore, as discussed in detail in Appendix E.2 of \cite{PP3}, any
condition sufficient for size controllability under Gaussianity of the
disturbance vector $\mathbf{U}$ also implies size controllability for large
classes of distributions for $\mathbf{U}$ that satisfy appropriate
domination conditions; however, the corresponding size-controlling critical
values may then differ from the size-controlling critical values that apply
under Gaussianity.

\section{Conclusion}

In the case of testing a \emph{single} restriction, we have shown that the
sufficient condition for size controllability of heteroskedasticity robust
test statistics in \cite{PP21} can be replaced by a weaker sufficient
condition that is also necessary. This allows one -- in the case of testing
a single restriction -- to resolve the question of existence of (finite)
size-controlling critical values in all cases, including those that remain
inconclusive under the results in \cite{PP21}.

We finally remark that the algorithms designed to compute size-controlling
critical values as discussed in Section 10 and Appendix E of \cite{PP21} can
be used as they stand also in situations where (a single restriction is
tested and) size controllability has been verified through checking
condition (\ref{non-incl_Het_uncorr}) (or equivalently (\ref{sol})) and
appealing to Theorem~\ref{thm:pp24}, but where (\ref{non-incl_Het}) does not
hold. This is so since the discussion of the before mentioned algorithms in 
\cite{PP21} only requires existence of a (finite) size-controlling critical
value, but does not depend on the way this existence is verified.

\appendix

\section{Auxiliary results and some examples \label{app:proofs}}

As a point of interest we note that Lemmata \ref{lem:3aux}, \ref%
{lem:additionV}, and \ref{gensol} below do \emph{not} rely on Assumption \ref%
{R_and_X}. Furthermore, all the lemmata in this appendix do neither refer to
the heteroskedasticity model nor to the Gaussianity assumption at all.
Finally, recall from Section \ref{Sec_2.2} that the set $\mathsf{B}$ is a
linear space (as $R$ is $1\times k$ in the present article).

\begin{lemma}
\label{lem:3aux} The following statements hold:.

\begin{enumerate}
\item $\mathsf{B}=\limfunc{span}(X)\oplus \{\hat{u}(y):y\in \mathsf{B}\}$,
the sum being orthogonal.

\item $\{\hat{u}(y):y\in \mathsf{B}\}$ is a linear subspace of $\limfunc{span%
}(e_{i}(n):i\in \mathcal{I}_{\#})$.

\item For every $z\in \limfunc{span}(e_{i}(n):i\in \mathcal{I}_{\#})$ we
have $R\hat{\beta}(z)=0$.

\item If $j\in \mathcal{I}_{\#}^{c}$, then $e_{j}(n)\in \limfunc{span}(X)$
and $e_{j}(n)\in \mathsf{B}$ are equivalent.
\end{enumerate}
\end{lemma}

\textbf{Proof:} 1. Obviously, $\{\hat{u}(y):y\in \mathsf{B}\}$ is a linear
space, since $\mathsf{B}$ is so. Observe that $\hat{u}(\hat{u}(y))=\hat{u}%
(y) $ holds, from which it follows that $B(y)=B(\hat{u}(y))$. Consequently, $%
y\in \mathsf{B}$ implies $\hat{u}(y)\in \mathsf{B}$. Since $\mathsf{B}$ is
invariant under addition of elements of $\limfunc{span}(X)$, we obtain $%
\mathsf{B}\supseteq \limfunc{span}(X)\oplus \{\hat{u}(y):y\in \mathsf{B}\}$,
the sum obviously being orthogonal. For the reverse inclusion, write $y\in 
\mathsf{B}$ as $y=X\hat{\beta}(y)+\hat{u}(y)$, which immediately implies
that $y\in \limfunc{span}(X)\oplus \{\hat{u}(y):y\in \mathsf{B}\}$.

2. Let $y\in \mathsf{B}$, i.e., $B(y)=0$, or, in other words, $R(X^{\prime
}X)^{-1}x_{i\cdot }^{\prime }\hat{u}_{i}(y)=0$ for every $i=1,...,n$. It
follows that $\hat{u}_{i}(y)=0$ for every $i\notin \mathcal{I}_{\#}$, from
which we conclude $\hat{u}(y)\in \limfunc{span}(e_{i}(n):i\in \mathcal{I}%
_{\#})$.

3. With $z_{i}$ denoting the $i$-th coordinate of $z$, we have%
\begin{eqnarray*}
R\hat{\beta}(z) &=&R(X^{\prime }X)^{-1}X^{\prime }z=R(X^{\prime
}X)^{-1}\dsum\limits_{i=1}^{n}z_{i}x_{i\cdot }^{\prime
}=\dsum\limits_{i=1}^{n}z_{i}R(X^{\prime }X)^{-1}x_{i\cdot }^{\prime } \\
&=&\dsum\limits_{i\in \mathcal{I}_{\#}}z_{i}R(X^{\prime }X)^{-1}x_{i\cdot
}^{\prime }+\dsum\limits_{i\in \mathcal{I}_{\#}^{c}}z_{i}R(X^{\prime
}X)^{-1}x_{i\cdot }^{\prime }=0,
\end{eqnarray*}%
observing that $R(X^{\prime }X)^{-1}x_{i\cdot }^{\prime }=0$ for $i\in 
\mathcal{I}_{\#}$ and that $z_{i}=0$ for $i\in \mathcal{I}_{\#}^{c}$.

4. Follows from the first two claims upon noting that~$j\in \mathcal{I}%
_{\#}^{c}$ is equivalent to $e_{j}(n)\bot \limfunc{span}(e_{i}(n):i\in 
\mathcal{I}_{\#})$. $\blacksquare $

\bigskip

\textbf{Remark A.1:} We discuss a few simple consequences of the preceding
lemma.

(i) If $\mathcal{I}_{\#}$ is empty then $\mathsf{B}=\limfunc{span}(X)$.

(ii) If $\mathcal{I}_{\#}=\{i_{0}\}$, then $\mathsf{B}=\limfunc{span}(X)$ or 
$\mathsf{B}=\limfunc{span}(X)\oplus \limfunc{span}(e_{i_{0}}(n))$; the
former happens if the $i_{0}$-th row of $X$ is nonzero, and the latter
happens if this row is zero.

(iii) If $\mathcal{I}_{\#}$ contains more than one element, then $\mathsf{B}=%
\limfunc{span}(X)$ (see (iv) below) as well as $\mathsf{B}\supsetneq 
\limfunc{span}(X)$ (see Example A.1 below) can occur.

(iv) Suppose $k=n-1$ and that Assumption \ref{R_and_X} holds. Then $\mathsf{B%
}=\limfunc{span}(X)$ always holds (since $\mathsf{B}$ is a linear space
containing the $n-1$ dimensional subspace $\limfunc{span}(X)$ and since $%
\mathsf{B}$ must be a proper subspace under Assumption \ref{R_and_X}, see
Lemma 3.1 in \cite{PP21}) regardless of whether $\mathcal{I}_{\#}$ is empty
or not. [That $\mathcal{I}_{\#}$ can indeed be nonempty in this situation is
shown by the example where $n=4$, $k=3$, $R=(1,1,0)^{\prime }$, and $X$ has
columns $(1,1,1,1)^{\prime }$, $(1,-1,1,-1)^{\prime }$, and $%
(1,1,-1,-1)^{\prime }$. It is easy to see that $e_{i}(4)\notin \limfunc{span}%
(X)$ for every $i=1,\ldots ,4$, and thus Assumption \ref{R_and_X} is
satisfied. The set $\mathcal{I}_{\#}$ is easily computed to be $\{2,4\}$.]

\begin{lemma}
\label{lem:additionV} The following statements hold:

\begin{enumerate}
\item The map $B$ and the set $\mathsf{B}$ are invariant w.r.t.~addition of
elements of $\mathsf{B}$. In particular, they are invariant w.r.t.~addition
of elements of $\mathcal{L}_{\#}:=\limfunc{span}(\mathfrak{M}_{0}^{lin}\cup 
\mathcal{V}_{\#})$.

\item $T_{Het}$ is invariant w.r.t.~addition of any $z\in \mathsf{B}$ that
satisfies $R\hat{\beta}(z)=0$.

\item $T_{Het}$ is invariant w.r.t.~addition of elements of $\mathcal{L}%
_{\#}=\limfunc{span}(\mathfrak{M}_{0}^{lin}\cup \mathcal{V}_{\#})$.
\end{enumerate}
\end{lemma}

\textbf{Proof:} 1. Linearity of~$B:\mathbb{R}^{n}\rightarrow (\mathbb{R}%
^{n})^{\prime }$ together with~$B(z)=0$ for every $z\in \mathsf{B}$ proves
the first statement in Part 1. [The invariance claim regarding $\mathsf{B}$
also trivially follows since $\mathsf{B}$ is a linear space.] The second one
then follows since, noting that $\mathsf{B}$ being a linear space, $%
\mathfrak{M}_{0}^{lin}\subseteq \limfunc{span}(X)\subseteq $ $\mathsf{B}$
and (\ref{eqn:Vdef}) imply $\mathcal{L}_{\#}\subseteq \mathsf{B}$.

2. First note that for $y\in \mathbb{R}^{n}$ and $z\in \mathsf{B}$ we have $%
\hat{\Omega}_{Het}(y+z)=\hat{\Omega}_{Het}(y)$ which follows from the easily
checked representation $\hat{\Omega}_{Het}(\cdot )=B(\cdot )\limfunc{diag}%
(d_{1},\ldots ,d_{n})B^{\prime }(\cdot )$ and Part 1 of the present lemma.
Second, clearly $R\hat{\beta}(y+z)-r=R\hat{\beta}(y)+R\hat{\beta}(z)-r=R\hat{%
\beta}(y)-r$ holds for $z$ satisfying $R\hat{\beta}(z)=0$. The claim now
follows from the definition of $T_{Het}$.

3. Follows from Part 2, since $\mathcal{L}_{\#}$ is a subset of $\mathsf{B}$
as shown in the proof of Part 1 of the present lemma, and since $z\in 
\mathcal{L}_{\#}$ implies $R\hat{\beta}(z)=0$ (because of linearity of $R%
\hat{\beta}(\cdot )$, because of the definition of $\mathfrak{M}_{0}^{lin}$,
and because of $\mathcal{V}_{\#}\subseteq \limfunc{span}(e_{i}(n):i\in 
\mathcal{I}_{\#})$ together with Part 3 of Lemma \ref{lem:3aux}). $%
\blacksquare $

\begin{lemma}
\label{lem:dimbound} Under Assumption \ref{R_and_X} we have $\limfunc{dim}(%
\mathcal{L}_{\#})<n-1.$
\end{lemma}

\textbf{Proof:} As shown in the proof of Part 1 of Lemma \ref{lem:additionV}%
, the relation $\mathcal{L}_{\#}\subseteq \mathsf{B}$ holds. Because $%
\mathsf{B}$ is a proper linear subspace of $\mathbb{R}^{n}$ under Assumption %
\ref{R_and_X} (cf.~Lemma 3.1 in \cite{PP21} and note that we have $q=1$
here), we must have $\limfunc{dim}(\mathcal{L}_{\#})\leq n-1$.\footnote{%
Alternatively, $\limfunc{dim}(\mathcal{L}_{\#})=n$ and invariance under
addition of elements of $\mathcal{L}_{\#}$ would lead to constancy of $%
T_{Het}$, and thus to a contradicition similar to the one arrived at in the
proof in the case $\limfunc{dim}(\mathcal{L}_{\#})=n-1$.} Assume now that $%
\mathcal{L}_{\#}$ has dimension $n-1$. Denote by $v\neq 0$ a vector that
spans $\mathcal{L}_{\#}^{\bot }$, the orthogonal complement of $\mathcal{L}%
_{\#}$ in $\mathbb{R}^{n}$, and fix an arbitrary $\mu _{0}\in \mathfrak{M}%
_{0}$. Use the invariance property in Part 3 of Lemma \ref{lem:additionV} to
see that for every $y\notin \mathcal{L}_{\#}$ we can write 
\begin{equation*}
T_{Het}(\mu _{0}+y)=T_{Het}(\mu _{0}+\Pi _{\mathcal{L}_{\#}^{\bot
}}y)=T_{Het}(\mu _{0}+v),
\end{equation*}%
where we used $\Pi _{\mathcal{L}_{\#}^{\bot }}y\neq 0$ together with
invariance of $T_{Het}$ w.r.t.~$G(\mathfrak{M}_{0})$ (cf.~Remark C.1 in
Appendix C of \cite{PP21}) to conclude the second equality.\footnote{%
Since $y\notin \mathcal{L}_{\#}$ we have $\Pi _{\mathcal{L}_{\#}^{\bot
}}y\neq 0,$ and thus $\Pi _{\mathcal{L}_{\#}^{\bot }}y=\lambda v$ with $%
\lambda \neq 0$. Invariance w.r.t. the group $G(\mathfrak{M}_{0})$ then
gives $T_{Het}(\mu _{0}+v)=T_{Het}(\mu _{0}+\lambda v)$.} But this implies
that $T_{Het}(\cdot )=T_{Het}(\mu _{0}+v)$ almost everywhere w.r.t.~Lebesgue
measure on $\mathbb{R}^{n}$, contradicting Part 2 of Lemma 5.16 in \cite{PP3}
in view of Remark C.1 in Appendix C of \cite{PP21} and noting that
Assumption \ref{R_and_X} is being maintained.\footnote{%
That $\limfunc{dim}(\mathcal{L}_{\#})=n-1$ leads to Lebesgue almost
everywhere constancy has been noted in Remark 5.14(i) of \cite{PP3} for a
large class of test statistics. We have included a proof here for the
convenience of the reader.} $\blacksquare $

\bigskip

\textbf{Remark A.2:} Even without Assumption \ref{R_and_X} we always have $%
\limfunc{dim}(\mathcal{L}_{\#})<n$. To see this, note that $I_{0}(\mathcal{L}%
_{\#})$ is a proper subset of $\{1,\ldots ,n\}$ by Part 3 of Lemma \ref%
{gensol} below, and thus $I_{1}(\mathcal{L}_{\#})\neq \emptyset $. But this
means that $e_{i}(n)\notin \mathcal{L}_{\#}$ for at least one $i$,
establishing the claim.

\begin{lemma}
\label{gensol} The following statements hold:

\begin{enumerate}
\item $i\in \mathcal{I}_{\#}$ if and only if $\Pi _{\limfunc{span}%
(X)}e_{i}(n)\in \mathfrak{M}_{0}^{lin}$.

\item Suppose $e_{i}(n)\in \limfunc{span}(X)$. Then $i\in \mathcal{I}_{\#}$
if and only if $i\in I_{0}(\mathfrak{M}_{0}^{lin})$.

\item $I_{0}(\mathfrak{M}_{0}^{lin})\subseteq I_{0}(\mathcal{L}%
_{\#})\subseteq \mathcal{I}_{\#}$ holds, and $\mathcal{I}_{\#}$ is a proper
subset of $\{1,\ldots ,n\}$.
\end{enumerate}
\end{lemma}

\textbf{Proof:} 1. Observe that%
\begin{eqnarray*}
R(X^{\prime }X)^{-1}x_{i\cdot }^{\prime } &=&R(X^{\prime }X)^{-1}X^{\prime
}e_{i}(n)=R(X^{\prime }X)^{-1}X^{\prime }(\Pi _{\limfunc{span}%
(X)}e_{i}(n)+\Pi _{\limfunc{span}(X)^{\bot }}e_{i}(n)) \\
&=&R(X^{\prime }X)^{-1}X^{\prime }\Pi _{\limfunc{span}(X)}e_{i}(n)=R\gamma
^{(i)},
\end{eqnarray*}%
where $\gamma ^{(i)}\in \mathbb{R}^{k}$ satisfies $\Pi _{\limfunc{span}%
(X)}e_{i}(n)=X\gamma ^{(i)}$. Consequently, $i\in \mathcal{I}_{\#}$ (i.e., $%
R(X^{\prime }X)^{-1}x_{i\cdot }^{\prime }=0$) if and only if $R\gamma
^{(i)}=0$ which is tantamount to $\Pi _{\limfunc{span}(X)}e_{i}(n)\in 
\mathfrak{M}_{0}^{lin}$.

2. Follows immediately from Part 1 and the definition of $I_{0}(\mathfrak{M}%
_{0}^{lin})$ upon noting that $\Pi _{\limfunc{span}(X)}e_{i}(n)=e_{i}(n)$
because of the assumption $e_{i}(n)\in \limfunc{span}(X)$.

3. The first inclusion is trivial since $\mathfrak{M}_{0}^{lin}\subseteq 
\mathcal{L}_{\#}$. To prove the second inclusion, suppose $i\in I_{0}(%
\mathcal{L}_{\#})$. Then $e_{i}(n)\in \mathcal{L}_{\#}$, which implies that $%
e_{i}(n)=v+w$ where $v\in \mathcal{V}_{\#}$ and $w\in \mathfrak{M}_{0}^{lin}$
(here we also use that $\mathcal{V}_{\#}$ and $\mathfrak{M}_{0}^{lin}$ are
linear subspaces). Using the definition of $\mathcal{V}_{\#}$ we arrive at%
\begin{equation*}
e_{i}(n)=\sum_{j:j\in \mathcal{I}_{\#},e_{j}(n)\in \mathsf{B}}\lambda
_{j}e_{j}(n)+w.
\end{equation*}%
Taking the projection and noting that $\Pi _{\limfunc{span}(X)}w=w$ (since $%
w\in \mathfrak{M}_{0}^{lin}\subseteq \limfunc{span}(X)$) this gives%
\begin{equation*}
\Pi _{\limfunc{span}(X)}e_{i}(n)=\sum_{j:j\in \mathcal{I}_{\#},e_{j}(n)\in 
\mathsf{B}}\lambda _{j}\Pi _{\limfunc{span}(X)}e_{j}(n)+w.
\end{equation*}%
The already established Part 1 shows that $\Pi _{\limfunc{span}%
(X)}e_{j}(n)\in \mathfrak{M}_{0}^{lin}$ for $j\in \mathcal{I}_{\#}$. Since $%
\mathfrak{M}_{0}^{lin}$ is a linear space we conclude that $\Pi _{\limfunc{%
span}(X)}e_{i}(n)$ belongs to $\mathfrak{M}_{0}^{lin}$. Again using Part 1,
we arrive at $i\in \mathcal{I}_{\#}$. That $\mathcal{I}_{\#}$ is a proper
subset of $\{1,\ldots ,n\}$ follows since $R\neq 0$. $\blacksquare $

\bigskip

\textbf{Remark A.3:} (i) Example A.1 below and the example discussed towards
the end of Remark A.1(iv) show that the first two inclusions in Part 3 of
the above lemma can be strict inclusions.

(ii) Inspection of the proof shows that Lemma \ref{gensol} actually also
holds if, in the notation of \cite{PP21}, we have $q\geq 1$, i.e., if a
collection of $q$ restrictions is tested simultaneously.

\bigskip

The subsequent examples show that condition (\ref{non-incl_Het}) can be
stronger than condition (\ref{non-incl_Het_uncorr}), another such example
being Example C.1 in Appendix C.1 of \cite{PP21}. We provide four different
examples to illustrate that this can happen in a variety of different
situations (e.g., independently of whether standard basis vectors belong to $%
\limfunc{span}(X)$ or not, etc.). We also compute the set $\mathsf{B}$ in
the examples below and illustrate the results in Lemma \ref{lem:3aux}.

\bigskip

\textbf{Example A.1:} Suppose $k=2$, $n=4$, and $X$ has $(1,1,1,1)^{\prime }$
as its first column and $(1,-1,1,-1)^{\prime }$ as its second column. Define
the $1\times k$ vector $R=(1,1)$. Then $\limfunc{rank}(X)=k=2$ holds, and $%
e_{j}(4)\notin \limfunc{span}(X)$ for every $j=1,\ldots ,4$, as is easily
checked; in particular, Assumption \ref{R_and_X} is thus satisfied, and $%
I_{1}(\mathfrak{M}_{0}^{lin})=\{1,\ldots ,4\}$. Furthermore, $R(X^{\prime
}X)^{-1}x_{i\cdot }^{\prime }\neq 0$ for $i=1,3$ whereas $R(X^{\prime
}X)^{-1}x_{i\cdot }^{\prime }=0$ for $i=2,4$. I.e., $\mathcal{I}%
_{\#}=\{2,4\} $. Now, $y\in \mathsf{B}$ (i.e., $B(y)=0$) is easily seen to
be equivalent to $\hat{u}_{1}(y)=\hat{u}_{3}(y)=0$, which in turn is
equivalent to $y_{1}=y_{3}$. In particular, $e_{2}(4)$ and $e_{4}(4)$ belong
to $\mathsf{B}$, but do not belong to $\limfunc{span}(X)$, while $e_{1}(4)$
and $e_{3}(4)$ do not belong to $\mathsf{B}$. The space $\{\hat{u}(y):y\in 
\mathsf{B}\}$ in the orthogonal sum representation $\mathsf{B}=\limfunc{span}%
(X)\oplus \{\hat{u}(y):y\in \mathsf{B}\}$ is here given by $\limfunc{span}%
((0,1,0,-1)^{\prime })$ as is not difficult to see. Note that, while $%
e_{2}(4)$ and $e_{4}(4)$ belong to $\mathsf{B}$ (and trivially also to $%
\limfunc{span}(e_{i}(4):i\in \mathcal{I}_{\#})$), they are not orthogonal to 
$\limfunc{span}(X)$, and do not belong to $\limfunc{span}((0,1,0,-1)^{\prime
})$ (which is a subset of $\limfunc{span}(e_{i}(4):i\in \mathcal{I}_{\#})$).
Furthermore, since $I_{1}(\mathfrak{M}_{0}^{lin})=\{1,\ldots ,4\}$,
condition (\ref{non-incl_Het_uncorr}) is satisfied, while condition (\ref%
{non-incl_Het}) is not. Theorem \ref{thm:pp21} does not allow one to draw a
conclusion about size-controllability of $T_{Het}$ in this example, while
Theorem \ref{thm:pp24} shows that $T_{Het}$ is size-controllable.

\textbf{Example A.2:} Suppose $k=3$, $n=5$, and $X$ has $(1,1,1,1,0)^{\prime
}$ as its first column, $(1,-1,1,-1,0)^{\prime }$ as its second column, and $%
(0,0,0,0,2)^{\prime }$ as its last column. Define the $1\times k$ vector $%
R=(1,1,r_{3})$. Then $\limfunc{rank}(X)=k=3$ holds, and $e_{j}(5)\notin 
\limfunc{span}(X)$ for every $j=1,\ldots ,4$, but $e_{5}(5)\in \limfunc{span}%
(X)$. Assumption \ref{R_and_X} is satisfied as can be easily checked.
Furthermore, $R(X^{\prime }X)^{-1}x_{i\cdot }^{\prime }\neq 0$ for $i=1,3$,
whereas $R(X^{\prime }X)^{-1}x_{i\cdot }^{\prime }=0$ for $i=2,4$; and $%
R(X^{\prime }X)^{-1}x_{5\cdot }^{\prime }=r_{3}/2$. Hence, $\mathcal{I}%
_{\#}=\{2,4\}$ in case $r_{3}\neq 0$, and $\mathcal{I}_{\#}=\{2,4,5\}$
otherwise. Now, $y\in \mathsf{B}$ (i.e., $B(y)=0$) is easily seen to be
equivalent to $\hat{u}_{1}(y)=\hat{u}_{3}(y)=0$, which in turn is equivalent
to $y_{1}=y_{3}$. In particular, $e_{2}(5)$ and $e_{4}(5)$ belong to $%
\mathsf{B}$, but do not belong to $\limfunc{span}(X)$, while $e_{5}(5)\in 
\limfunc{span}(X)\subseteq \mathsf{B}$; and $e_{1}(5)$ and $e_{3}(5)$ do not
belong to $\mathsf{B}$. The space $\{\hat{u}(y):y\in \mathsf{B}\}$ in the
orthogonal sum representation $\mathsf{B}=\limfunc{span}(X)\oplus \{\hat{u}%
(y):y\in \mathsf{B}\}$ is here given by $\limfunc{span}((0,1,0,-1,0)^{\prime
})$ as is not difficult to see. Note that, while $e_{2}(5)$ and $e_{4}(5)$
belong to $\mathsf{B}$ (and trivially also to $\limfunc{span}(e_{i}(5):i\in 
\mathcal{I}_{\#})$), they are not orthogonal to $\limfunc{span}(X)$, and do
not belong to $\limfunc{span}((0,1,0,-1,0)^{\prime })$ (which is a subset of 
$\limfunc{span}(e_{i}(5):i\in \mathcal{I}_{\#})$). Note that $I_{1}(%
\mathfrak{M}_{0}^{lin})=\{1,\ldots ,4\}$ in case $r_{3}=0$, while $I_{1}(%
\mathfrak{M}_{0}^{lin})=\{1,\ldots ,5\}$ otherwise. In particular, in case $%
r_{3}=0$, condition (\ref{non-incl_Het_uncorr}) is satisfied, while
condition (\ref{non-incl_Het}) is not; hence, in this case Theorem \ref%
{thm:pp21} does not allow one to draw a conclusion about
size-controllability of $T_{Het}$, while Theorem \ref{thm:pp24} shows that $%
T_{Het}$ is size-controllable. In case $r_{3}\neq 0$, both conditions (\ref%
{non-incl_Het}) and (\ref{non-incl_Het_uncorr}) are violated, and both
theorems show that the test based on $T_{Het}$ has size $1$ regardless of
the choice of critical value.

\textbf{Example A.3:} Suppose $k=2$, $n=5$, and $X$ has $(1,1,1,1,0)^{\prime
}$ as its first column and $(1,-1,1,-1,0)^{\prime }$ as its second column.
Define the $1\times k$ vector $R=(1,0)$. Then $\limfunc{rank}(X)=k=2$ holds,
and $e_{j}(5)\notin \limfunc{span}(X)$ for every $j=1,\ldots ,5$, as is
easily checked; in particular, Assumption \ref{R_and_X} is thus satisfied,
and $I_{1}(\mathfrak{M}_{0}^{lin})=\{1,\ldots ,5\}$. Furthermore, $%
R(X^{\prime }X)^{-1}x_{i\cdot }^{\prime }\neq 0$ for $i=1,\ldots ,4$ whereas 
$R(X^{\prime }X)^{-1}x_{5\cdot }^{\prime }=0$. I.e., $\mathcal{I}_{\#}=\{5\}$%
. Now, $y\in \mathsf{B}$ (i.e., $B(y)=0$) is easily seen to be equivalent to 
$\hat{u}_{1}(y)=\hat{u}_{2}(y)=\hat{u}_{3}(y)=\hat{u}_{4}(y)=0$, which in
turn is equivalent to $y_{1}=y_{3}$ and $y_{2}=y_{4}$. In particular, $%
e_{5}(5)$ belongs to $\mathsf{B}$, but does not belong to $\limfunc{span}(X)$%
, in fact is orthogonal to $\limfunc{span}(X)$, while $e_{j}(5)\notin 
\mathsf{B}$ for $j=1,\ldots ,4$. The space $\{\hat{u}(y):y\in \mathsf{B}\}$
in the orthogonal sum representation $\mathsf{B}=\limfunc{span}(X)\oplus \{%
\hat{u}(y):y\in \mathsf{B}\}$ is here given by $\limfunc{span}(e_{5}(5))$ as
is not difficult to see. Furthermore, in this example condition (\ref%
{non-incl_Het_uncorr}) is satisfied, while condition (\ref{non-incl_Het}) is
not. Theorem \ref{thm:pp21} does not allow one to draw a conclusion about
size-controllability of $T_{Het}$ in this example, while Theorem \ref%
{thm:pp24} shows that $T_{Het}$ is size-controllable.

\textbf{Example A.4:} Suppose $k=3$, $n=6$, and $X$ has $(1,1,1,1,0,0)^{%
\prime }$ as its first column, $(1,-1,1,-1,0,0)^{\prime }$ as its second
column, and $(0,0,0,0,0,2)^{\prime }$ as its third column. Define the $%
1\times k$ vector $R=(1,0,r_{3})$. Then $\limfunc{rank}(X)=k=3$ holds, and $%
e_{j}(6)\notin \limfunc{span}(X)$ for every $j=1,\ldots ,5$, but $%
e_{6}(6)\in \limfunc{span}(X)$. Assumption \ref{R_and_X} is satisfied as can
be easily checked. Furthermore, $R(X^{\prime }X)^{-1}x_{i\cdot }^{\prime
}\neq 0$ for $i=1,\ldots ,4$ whereas $R(X^{\prime }X)^{-1}x_{5\cdot
}^{\prime }=0$ and $R(X^{\prime }X)^{-1}x_{6\cdot }^{\prime }=r_{3}/2$.
Hence, $\mathcal{I}_{\#}=\{5\}$ in case $r_{3}\neq 0,$ and $\mathcal{I}%
_{\#}=\{5,6\}$ otherwise. Now, $y\in \mathsf{B}$ (i.e., $B(y)=0$) is easily
seen to be equivalent to $\hat{u}_{1}(y)=\hat{u}_{2}(y)=\hat{u}_{3}(y)=\hat{u%
}_{4}(y)=0$, which in turn is equivalent to $y_{1}=y_{3}$ and $y_{2}=y_{4}$.
In particular, $e_{5}(6)$ belongs to $\mathsf{B}$, but does not belong to $%
\limfunc{span}(X)$, in fact is orthogonal to $\limfunc{span}(X)$, while $%
e_{6}(6)\in \limfunc{span}(X)\subseteq \mathsf{B}$; and $e_{j}(6)\notin 
\mathsf{B}$ for $j=1,\ldots ,4$. The space $\{\hat{u}(y):y\in \mathsf{B}\}$
in the orthogonal sum representation $\mathsf{B}=\limfunc{span}(X)\oplus \{%
\hat{u}(y):y\in \mathsf{B}\}$ is here given by $\limfunc{span}(e_{5}(6))$ as
is not difficult to see. Note that $I_{1}(\mathfrak{M}_{0}^{lin})=\{1,\ldots
,5\}$ in case $r_{3}=0$, while $I_{1}(\mathfrak{M}_{0}^{lin})=\{1,\ldots
,6\} $ otherwise. In particular, in case $r_{3}=0$, condition (\ref%
{non-incl_Het_uncorr}) is satisfied, while condition (\ref{non-incl_Het}) is
not; hence, in this case Theorem \ref{thm:pp21} does not allow one to draw a
conclusion about size-controllability of $T_{Het}$, while Theorem \ref%
{thm:pp24} shows that $T_{Het}$ is size-controllable. In case $r_{3}\neq 0$,
both conditions (\ref{non-incl_Het}) and (\ref{non-incl_Het_uncorr}) are
violated, and both theorems show that the test based on $T_{Het}$ has size $%
1 $ regardless of the choice of critical value.

\bigskip

\textbf{Remark A.4:} Many more examples can be generated from Examples
A.1-A.4 via the transformation $X^{\ast }=XA$ and $R^{\ast }=RA$ where $A$
is a nonsingular $k\times k$ matrix. These new examples exhibit the same
features as Examples A.1-A.4, respectively. In particular, one can generate
examples that have $R^{\ast }=(1,0\ldots ,0)$.

\section{Proof of Theorem \protect\ref{thm:pp24}\label{app:B}}

To prove Theorem \ref{thm:pp24} we follow the strategy used to establish
Theorem 5.1 in \cite{PP21} and first provide a result for a class of
heteroskedasticity models that includes $\mathfrak{C}_{Het}$ as a special
case, and which is of some independent interest. The heteroskedasticity
models we consider here are defined as follows (cf.~Appendix A of \cite{PP21}
for more discussion): Let $m\in \mathbb{N}$, and let $n_{j}\in \mathbb{N}$
for $j=1,\ldots ,m$ satisfy $\sum_{j=1}^{m}n_{j}=n$. Set $%
n_{j}^{+}=\sum_{l=1}^{j}n_{l}$ and define%
\begin{equation}
\mathfrak{C}_{(n_{1},\ldots ,n_{m})}=\left\{ \limfunc{diag}(\tau
_{1}^{2},\ldots ,\tau _{n}^{2})\in \mathfrak{C}_{Het}:\tau
_{n_{j-1}^{+}+1}^{2}=\ldots =\tau _{n_{j}^{+}}^{2}\text{ for }j=1,\ldots
,m\right\}  \label{eqn:covmodmg}
\end{equation}%
with the convention that $n_{0}^{+}=0$. In the special case where $m=n$ and $%
n_{1}=n_{2}=...=n_{m}=1$ we have $\mathfrak{C}_{(n_{1},\ldots ,n_{m})}=%
\mathfrak{C}_{Het}$. We use $\lambda _{\mathbb{R}^{n}}$ to denote Lebesgue
measure on $\mathbb{R}^{n}$, and $\lambda _{\mathcal{A}}$ to denote Lebesgue
measure on a (nonempty) affine space $\mathcal{A}$ (but viewed as a measure
on the Borel-sets of $\mathbb{R}^{n}$), with zero-dimensional Lebesgue
measure interpreted as point mass. We start with a lemma and note that it
does \emph{not} make use of Assumption \ref{R_and_X}. Recall that by
definition $\mathcal{L}_{\#}=\limfunc{span}(\mathfrak{M}_{0}^{lin}\cup 
\mathcal{V}_{\#})$, and that \emph{we only consider testing a single
restriction in the present article}.

\begin{lemma}
\label{equivalence} Let $m\in \mathbb{N}$, and let $n_{j}\in \mathbb{N}$ for 
$j=1,\ldots ,m$ satisfy $\sum_{j=1}^{m}n_{j}=n$. Then:

(a) The condition%
\begin{align}
\func{span}\left( \left\{ e_{i}(n):i\in (n_{j-1}^{+},n_{j}^{+}]\right\}
\right) & \nsubseteqq \mathsf{B}\text{ }  \notag \\
\text{\ for every }j& =1,\ldots ,m\text{ with }(n_{j-1}^{+},n_{j}^{+}]\cap
I_{1}(\mathcal{L}_{\#})\neq \emptyset  \label{non-incl_groupHet}
\end{align}%
is equivalent to the condition%
\begin{align}
\limfunc{span}\left( \left\{ e_{i}(n):i\in (n_{j-1}^{+},n_{j}^{+}]\cap I_{1}(%
\mathcal{L}_{\#}\right\} \right) & \nsubseteqq \limfunc{span}(X)  \notag \\
\text{ for every }j& =1,\ldots ,m\text{ with }\emptyset \neq
(n_{j-1}^{+},n_{j}^{+}]\cap I_{1}(\mathcal{L}_{\#})\subseteq \mathcal{I}%
_{\#}^{c}.  \label{non-incl_groupHet_uncorr}
\end{align}%
[It is understood here, that condition (\ref{non-incl_groupHet_uncorr}) is
satisfied if no $j$ with $\emptyset \neq (n_{j-1}^{+},n_{j}^{+}]\cap I_{1}(%
\mathcal{L}_{\#})\subseteq \mathcal{I}_{\#}^{c}$ exists.]

(b) In the special case where $m=n$ and $n_{1}=n_{2}=...=n_{m}=1$, (\ref%
{non-incl_groupHet_uncorr}) (as well as (\ref{non-incl_groupHet})) is
equivalent to (\ref{non-incl_Het_uncorr}).
\end{lemma}

\textbf{Proof:} (a) Recall from the proof of Part 1 of Lemma \ref%
{lem:additionV} that $\mathcal{L}_{\#}=\limfunc{span}(\mathfrak{M}%
_{0}^{lin}\cup \mathcal{V}_{\#})\subseteqq \mathsf{B}$. Therefore, $%
e_{i}(n)\notin \mathsf{B}$ is possible only if $i\in I_{1}(\mathcal{L}_{\#})$%
. Hence, in view of invariance of $\mathsf{B}$ w.r.t. addition of elements
of $\mathsf{B}$ (Lemma \ref{lem:additionV}), the condition in (\ref%
{non-incl_groupHet}) is equivalent to 
\begin{eqnarray}
\func{span}\left( \left\{ e_{i}(n):i\in (n_{j-1}^{+},n_{j}^{+}]\cap I_{1}(%
\mathcal{L}_{\#})\right\} \right) &\nsubseteqq &\mathsf{B}  \notag \\
\text{\ for every }j &=&1,\ldots ,m\text{ with }(n_{j-1}^{+},n_{j}^{+}]\cap
I_{1}(\mathcal{L}_{\#})\neq \emptyset .  \label{non-incl_groupHetmod2a}
\end{eqnarray}%
For $i\in \mathcal{I}_{\#}$ the condition $e_{i}(n)\in \mathsf{B}$ implies $%
e_{i}(n)\in \mathcal{V}_{\#}\subseteq \mathcal{L}_{\#}$, so that $i\notin
I_{1}(\mathcal{L}_{\#})$. In other words, $i\in I_{1}(\mathcal{L}_{\#})\cap 
\mathcal{I}_{\#}$ implies $e_{i}(n)\notin \mathsf{B}$. This shows that for
any $j$ with the property that $(n_{j-1}^{+},n_{j}^{+}]\cap I_{1}(\mathcal{L}%
_{\#})$ contains an element $i\in \mathcal{I}_{\#}$, the non-inclusion
relation in (\ref{non-incl_groupHetmod2a}) is automatically satisfied.
Hence, (\ref{non-incl_groupHetmod2a}) is equivalent to 
\begin{eqnarray}
\func{span}\left( \left\{ e_{i}(n):i\in (n_{j-1}^{+},n_{j}^{+}]\cap I_{1}(%
\mathcal{L}_{\#})\right\} \right) &\nsubseteqq &\mathsf{B}
\label{non-incl_groupHetmod2} \\
\text{for every }j &=&1,\ldots ,m\text{ with }\emptyset \neq
(n_{j-1}^{+},n_{j}^{+}]\cap I_{1}(\mathcal{L}_{\#})\subseteq \mathcal{I}%
_{\#}^{c}  \notag
\end{eqnarray}%
with the understanding that this condition is satisfied if no $j$ with $%
\emptyset \neq (n_{j-1}^{+},n_{j}^{+}]\cap I_{1}(\mathcal{L}_{\#})\subseteq 
\mathcal{I}_{\#}^{c}$ exists. Since $\mathsf{B}$ as well as $\limfunc{span}%
(X)$ are a linear spaces, Part 4 of Lemma~\ref{lem:3aux} shows that (\ref%
{non-incl_groupHetmod2}) is equivalent to the statement in (\ref%
{non-incl_groupHet_uncorr}).

(b) In the special case considered here (\ref{non-incl_groupHet_uncorr})
simplifies to 
\begin{equation}
e_{i}(n)\notin \limfunc{span}(X)\quad \text{ \ for every }i\in I_{1}(%
\mathcal{L}_{\#})\cap \mathcal{I}_{\#}^{c}
\label{non-incl_groupHetmod2spez*}
\end{equation}%
with the understanding as in (\ref{non-incl_groupHet_uncorr}) that this
condition is satisfied if $I_{1}(\mathcal{L}_{\#})\cap \mathcal{I}_{\#}^{c}$
is empty. Since $I_{1}(\mathcal{L}_{\#})\cap \mathcal{I}_{\#}^{c}=\mathcal{I}%
_{\#}^{c}\neq \emptyset $ by Part 3 of Lemma \ref{gensol}, the index set in (%
\ref{non-incl_groupHetmod2spez*}) is actually nonempty, and furthermore (\ref%
{non-incl_groupHetmod2spez*}) is equivalent to%
\begin{equation}
e_{i}(n)\notin \limfunc{span}(X)\quad \text{ \ for every }i\in \mathcal{I}%
_{\#}^{c}.  \label{delta}
\end{equation}%
Because of $\mathcal{I}_{\#}^{c}\subseteq I_{1}(\mathfrak{M}_{0}^{lin})$
(Lemma \ref{gensol}), the statement in (\ref{delta}) is implied by that in (%
\ref{non-incl_Het_uncorr}). To show that (\ref{delta}) implies (\ref%
{non-incl_Het_uncorr}), suppose (\ref{non-incl_Het_uncorr}) is violated,
i.e., there exists an $i\in I_{1}(\mathfrak{M}_{0}^{lin})$ such that $%
e_{i}(n)\in \limfunc{span}(X)$. It then follows that $R\hat{\beta}%
(e_{i}(n))\neq 0$ must hold. Since $R\hat{\beta}(e_{i}(n))=R(X^{\prime
}X)^{-1}x_{i\cdot }^{\prime }$, we conclude $i\in \mathcal{I}_{\#}^{c}$.
Hence, also (\ref{delta}) must be violated, a contradiction. $\blacksquare $

\bigskip

Parts 1-2 of the following statement provide -- in the context of testing a
single restriction -- a version of Theorem A.1(b) and the corresponding part
of Theorem A.1(c) in \cite{PP21}, while Part 3 corresponds to the
generalization of Proposition 5.5(b) mentioned after Theorem A.1 in \cite%
{PP21}. Part 4 of the subsequent theorem is a version of Proposition A.2(b)
in \cite{PP21}, and together with Part 1 shows that under Assumption \ref%
{R_and_X} the condition in (\ref{non-incl_groupHet}), or equivalently (\ref%
{non-incl_groupHet_uncorr}), is \emph{necessary and sufficient} for the
existence of a (finite) critical value that controls the size of $T_{Het}$
over the heteroskedasticity model $\mathfrak{C}_{(n_{1},\ldots ,n_{m})}$
when testing%
\begin{equation*}
H_{0}:\mu \in \mathfrak{M}_{0},\ 0<\sigma ^{2}<\infty ,\ \Sigma \in 
\mathfrak{C}_{(n_{1},\ldots ,n_{m})}\quad \text{ vs. }\quad H_{1}:\mu \in 
\mathfrak{M}_{1},\ 0<\sigma ^{2}<\infty ,\ \Sigma \in \mathfrak{C}%
_{(n_{1},\ldots ,n_{m})}.
\end{equation*}

\begin{theorem}
\label{thm:moregen}Let $m\in \mathbb{N}$, let $n_{j}\in \mathbb{N}$ for $%
j=1,\ldots ,m$ satisfy $\sum_{j=1}^{m}n_{j}=n$, and suppose Assumption \ref%
{R_and_X} is satisfied. Then the following statements hold:

\begin{enumerate}
\item For every $0<\alpha <1$ there exists a real number $C(\alpha )$ such
that 
\begin{equation}
\sup_{\mu _{0}\in \mathfrak{M}_{0}}\sup_{0<\sigma ^{2}<\infty }\sup_{\Sigma
\in \mathfrak{C}_{(n_{1},\ldots ,n_{m})}}P_{\mu _{0},\sigma ^{2}\Sigma
}(T_{Het}\geq C(\alpha ))\leq \alpha  \label{size-control_groupHet}
\end{equation}%
holds, provided that (\ref{non-incl_groupHet}) (or equivalently (\ref%
{non-incl_groupHet_uncorr})) holds. Furthermore, under condition (\ref%
{non-incl_groupHet}) (or equivalently (\ref{non-incl_groupHet_uncorr})),
even equality can be achieved in (\ref{size-control_groupHet}) by a proper
choice of $C(\alpha )$, provided $\alpha \in (0,\alpha ^{\ast }]\cap (0,1)$
holds, where 
\begin{equation*}
\alpha ^{\ast }=\sup_{C\in (C^{\ast },\infty )}\sup_{\Sigma \in \mathfrak{C}%
_{(n_{1},\ldots ,n_{m})}}P_{\mu _{0},\Sigma }(T_{Het}\geq C)
\end{equation*}%
is positive and where $C^{\ast }$ is defined as in Lemma 5.11 of \cite{PP3}
with $\mathfrak{C}=\mathfrak{C}_{(n_{1},\ldots ,n_{m})}$, $T=T_{Het}$, $%
N^{\dag }=\mathsf{B}$, $\mathcal{L}=\mathcal{L}_{\#}$, and $q=1$ (with
neither $\alpha ^{\ast }$ nor $C^{\ast }$ depending on the choice of $\mu
_{0}\in \mathfrak{M}_{0}$).

\item Suppose (\ref{non-incl_groupHet}) (or equivalently (\ref%
{non-incl_groupHet_uncorr})) is satisfied. Then a smallest critical value,
denoted by $C_{\Diamond }(\alpha )$, satisfying (\ref{size-control_groupHet}%
) exists for every $0<\alpha <1$. And $C_{\Diamond }(\alpha )$ is also the
smallest among the critical values leading to equality in (\ref%
{size-control_groupHet}) whenever such critical values exist.\footnote{%
The dependence of $C_{\Diamond }(\alpha )$ on the heteroskedasticity model
is not shown in the notation, In particular, $C_{\Diamond }(\alpha )$ in the
current theorem is not necessarily the same as $C_{\Diamond }(\alpha )$ in
the other theorems.}

\item Suppose (\ref{non-incl_groupHet}) (or equivalently (\ref%
{non-incl_groupHet_uncorr})) is satisfied. Then any $C(\alpha )$ satisfying (%
\ref{size-control_groupHet}) necessarily has to satisfy $C(\alpha )\geq
C^{\ast }$. In fact, for any $C<C^{\ast }$ we have $\sup_{\Sigma \in 
\mathfrak{C}_{(n_{1},\ldots ,n_{m})}}P_{\mu _{0},\sigma ^{2}\Sigma
}(T_{Het}\geq C)=1$ for every $\mu _{0}\in \mathfrak{M}_{0}$ and every $%
\sigma ^{2}\in (0,\infty )$.

\item If (\ref{non-incl_groupHet}) (or equivalently (\ref%
{non-incl_groupHet_uncorr})) is violated, then $\sup_{\Sigma \in \mathfrak{C}%
_{(n_{1},\ldots ,n_{m})}}P_{\mu _{0},\sigma ^{2}\Sigma }(T_{Het}\geq C)=1$
for \emph{every} choice of critical value $C$, every $\mu _{0}\in \mathfrak{M%
}_{0}$, and every $\sigma ^{2}\in (0,\infty )$ (implying that size equals $1$
for every $C$).\footnote{%
Cf. Footnote \ref{FN}.}
\end{enumerate}
\end{theorem}

The following proof adapts the proof of Theorem A.1 in \cite{PP21}.

\textbf{Proof of Theorem \ref{thm:moregen}:} We first prove Part 1. We apply
Part A of Proposition 5.12 of \cite{PP3} with $\mathfrak{C}=\mathfrak{C}%
_{(n_{1},\ldots ,n_{m})}$, $T=T_{Het}$, $\mathcal{L}=\mathcal{L}_{\#}$, and $%
\mathcal{V}=\mathcal{V}_{\#}$ (and $q=1$). First, note that $\dim (\mathcal{L%
}_{\#})<n-1<n$ because of Lemma \ref{lem:dimbound}. Second, under Assumption %
\ref{R_and_X}, $T_{Het}$ is a non-sphericity corrected F-type test statistic
with $N^{\ast }=\mathsf{B}$, which is a closed $\lambda _{\mathbb{R}^{n}}$%
-null set (see Remarks 3.2 and C.1 as well as Lemma 3.1 in \cite{PP21}); in
particular, $T_{Het}$ as well as $\mathsf{B}$ are invariant w.r.t. the group 
$G(\mathfrak{M}_{0})$. Furthermore, $T_{Het}$ as well as $\mathsf{B}$ are
invariant w.r.t. addition of elements of $\mathcal{V}_{\#}$ by Lemma~\ref%
{lem:additionV}. Hence, the general assumptions on $T=T_{Het}$, on $N^{\dag
}=N^{\ast }=\mathsf{B}$, on $\mathcal{V}=\mathcal{V}_{\#}$, as well as on $%
\mathcal{L}=\mathcal{L}_{\#}$ in Proposition 5.12 of \cite{PP3} are
satisfied in view of Part 1 of Lemma 5.16 in the same reference.

Next, observe that condition (\ref{non-incl_groupHet}) is equivalent to%
\begin{equation*}
\func{span}\left( \left\{ \Pi _{\mathcal{L}_{\#}^{\bot }}e_{i}(n):i\in
(n_{j-1}^{+},n_{j}^{+}]\right\} \right) \nsubseteqq \mathsf{B}
\end{equation*}%
for every $j=1,\ldots ,m$, such that $(n_{j-1}^{+},n_{j}^{+}]\cap I_{1}(%
\mathcal{L}_{\#})\neq \emptyset $, since $\Pi _{\mathcal{L}_{\#}^{\bot
}}e_{i}(n)$ and $e_{i}(n)$ differ only by an element of $\mathcal{L}_{\#}$
and since $\mathsf{B}+\mathcal{L}_{\#}=\mathsf{B}$ (which follows from Part
1 of Lemma \ref{lem:additionV}). In view of Proposition B.2 in Appendix B of 
\cite{PP21}, this implies that any $\mathcal{S}\in \mathbb{J}(\mathcal{L}%
_{\#},\mathfrak{C}_{(n_{1},\ldots ,n_{m})})$ is not contained in $\mathsf{B}$%
, and thus not in $N^{\dag }$. Using $\mathfrak{M}_{0}\subseteq \limfunc{span%
}(X)$ and $\mathsf{B}+\limfunc{span}(X)=\mathsf{B}$ (by Lemma 3.1(e) in \cite%
{PP21}), it follows that $\mu _{0}+\mathcal{S}\nsubseteqq \mathsf{B}=N^{\dag
}$ for every $\mu _{0}\in \mathfrak{M}_{0}$. Since $\mu _{0}+\mathcal{S}$ is
a (nonempty) affine space and $N^{\dag }=\mathsf{B}$ is a linear space
(recall that $R$ is $1\times k$), we may conclude (cf.~Corollary 5.6 in \cite%
{PP3} and its proof) that $\lambda _{\mu _{0}+\mathcal{S}}(N^{\dag })=0$ for
every $\mathcal{S}\in \mathbb{J}(\mathcal{L}_{\#},\mathfrak{C}%
_{(n_{1},\ldots ,n_{m})})$ and every $\mu _{0}\in \mathfrak{M}_{0}$. This
completes the verification of the assumptions of Proposition 5.12 in \cite%
{PP3} that are not specific to Part A (or Part B) of this proposition.

We next verify the assumptions specific to Part A of this proposition:
Assumption (a) is satisfied (even for every $C\in \mathbb{R}$) as a
consequence of Part 2 of Lemma 5.16 in \cite{PP3} and of Remark C.1(i) in
Appendix C of \cite{PP21}. And Assumption (b) in Part A follows from Lemma
5.19 of \cite{PP3}, since $T_{Het}$ results as a special case of the test
statistics $T_{GQ}$ defined in Section 3.4 of \cite{PP3} upon choosing $%
\mathcal{W}_{n}^{\ast }=n^{-1}\limfunc{diag}(d_{1},\ldots ,d_{n})$. Part A
of Proposition 5.12 of \cite{PP3} now immediately delivers claim (\ref%
{size-control_groupHet}), since $C^{\ast }<\infty $ as noted in that
proposition. That $C^{\ast }$ and $\alpha ^{\ast }$ do not depend on the
choice of $\mu _{0}\in \mathfrak{M}_{0}$ is an immediate consequence of $G(%
\mathfrak{M}_{0})$-invariance of $T_{Het}$ (cf.~Remark 3.2 in \cite{PP21}).
Also note that $\alpha ^{\ast }$ as defined in the theorem coincides with $%
\alpha ^{\ast }$ as defined in Proposition 5.12 of \cite{PP3} in view of $G(%
\mathfrak{M}_{0})$-invariance of $T_{Het}$. Positivity of $\alpha ^{\ast }$
then follows from Part 5 of Lemma 5.15 in \cite{PP2016} in view of Remark
C.1(i) in Appendix C of \cite{PP21}, noting that $\lambda _{\mathbb{R}^{n}}$
and $P_{\mu _{0},\Sigma }$ are equivalent measures (since $\Sigma \in 
\mathfrak{C}_{Het}$ is positive definite); cf.~Remark 5.13(vi) in \cite{PP3}%
. In case $\alpha <\alpha ^{\ast }$, the remaining claim in Part 1 of the
present theorem, namely that equality can be achieved in (\ref%
{size-control_groupHet}), follows from the definition of $C^{\ast }$ in
Lemma 5.11 of \cite{PP3} and from Part A.2 of Proposition 5.12 of \cite{PP3}
(and the observation immediately following that proposition allowing one to
drop the suprema w.r.t. $\mu _{0}$ and $\sigma ^{2}$, and to set $\sigma
^{2}=1$); in case $\alpha =\alpha ^{\ast }<1$, it follows from Remarks
5.13(i),(ii) in \cite{PP3} using Lemma 5.16 in the same reference.

The claim in Part 2 follows from Remark 5.10 and Lemma 5.16 in \cite{PP3}
combined with Remark C.1(i) in Appendix C of \cite{PP21}; cf.~also Appendix
A.3 in \cite{PP21}.

Part 3 follows from Part A.1 of Proposition 5.12 of \cite{PP3} and the
sentence following this proposition. Note that the assumptions of this
proposition have been verified in the proof of Part 1 above.

Part 4 follows from Part 3 of Corollary 5.17 in \cite{PP2016}: As shown in
Remark C.1 in Appendix C of \cite{PP21}, $T_{Het}$ satisfies the assumptions
of this corollary (with $\check{\beta}=\hat{\beta}$, $\check{\Omega}=\hat{%
\Omega}_{Het}$, $N=\emptyset $, and $N^{\ast }=\mathsf{B}$). Suppose that (%
\ref{non-incl_groupHet_uncorr}) is violated and set $\mathcal{Z}=\func{span}%
(\{e_{i}(n):i\in (n_{j-1}^{+},n_{j}^{+}]\})$, where $j$ is such that $%
\emptyset \neq (n_{j-1}^{+},n_{j}^{+}]\cap I_{1}(\mathcal{L}_{\#})\subseteq 
\mathcal{I}_{\#}^{c}$ and 
\begin{equation}
\limfunc{span}\left( \left\{ e_{i}(n):i\in (n_{j-1}^{+},n_{j}^{+}]\cap I_{1}(%
\mathcal{L}_{\#})\right\} \right) \subseteq \limfunc{span}(X).  \label{cond}
\end{equation}%
Since $e_{i}(n)\in \mathcal{L}_{\#}$ for every $i\in I_{0}(\mathcal{L}_{\#})$%
, it hence follows from (\ref{cond}) that $\mathcal{Z}\subseteq \limfunc{span%
}(\limfunc{span}(X)\cup \mathcal{L}_{\#})\subseteq \mathsf{B}$, recalling
that $\limfunc{span}(X)\subseteq \mathsf{B}$, that $\mathcal{L}%
_{\#}\subseteq \mathsf{B}$ (cf.~the proof of Part 1 of Lemma \ref%
{lem:additionV}), and that $\mathsf{B}$ is a linear space (recall that $R$
is $1\times k$). Note that $\mathcal{Z}$ is not contained in $\mathfrak{M}%
_{0}^{lin}$ because $\emptyset \neq (n_{j-1}^{+},n_{j}^{+}]\cap I_{1}(%
\mathcal{L}_{\#})$ but $\mathfrak{M}_{0}^{lin}\subseteq \mathcal{L}_{\#}$.
Observe that $\mathcal{Z}$ is a concentration space of $\mathfrak{C}%
_{(n_{1},\ldots ,n_{m})}$ in view of Remark B.4 in Appendix B of \cite{PP21}
(note that $\limfunc{card}((n_{j-1}^{+},n_{j}^{+}])<n$ must hold in view of $%
\mathcal{Z}\subseteq \mathsf{B}$ and $\mathsf{B}$ being a proper subspace of 
$\mathbb{R}^{n}$ by Lemma 3.1 in \cite{PP21} in conjunction with Assumption %
\ref{R_and_X}, while $0<\limfunc{card}((n_{j-1}^{+},n_{j}^{+}])$ is
obvious). The nonnegative definiteness assumption on $\check{\Omega}=\hat{%
\Omega}_{Het}$ in Part 3 of Corollary 5.17 in \cite{PP2016} is satisfied
(cf.~Lemma 3.1 in \cite{PP21}). Obviously $\check{\Omega}(z)=0$ holds for
every $z\in \mathcal{Z}$ as a consequence of Part (b) of Lemma 3.1 in \cite%
{PP21} since $\mathcal{Z}\subseteq \mathsf{B}$ (as just shown) and since $%
\check{\Omega}(z)$ is $1\times 1$. It remains to establish that $R\hat{\beta}%
(z)\neq 0$ holds $\lambda _{\mathcal{Z}}$-everywhere: we recall that $%
\emptyset \neq (n_{j-1}^{+},n_{j}^{+}]\cap I_{1}(\mathcal{L}_{\#})\subseteq 
\mathcal{I}_{\#}^{c}$ and pick an element $i$, say, of $%
(n_{j-1}^{+},n_{j}^{+}]\cap I_{1}(\mathcal{L}_{\#})$. Then $e_{i}(n)\in 
\mathcal{Z}$ and $i\in \mathcal{I}_{\#}^{c}$, and from the definition of $%
\mathcal{I}_{\#}^{c}$ we conclude that $R\hat{\beta}(e_{i}(n))\neq 0$. It
follows that the linear space $\mathcal{Z}$ is not a subspace of the kernel
of $R\hat{\beta}$ so that $R\hat{\beta}(z)\neq 0$ holds $\lambda _{\mathcal{Z%
}}$-everywhere. Part 3 of Corollary 5.17 in \cite{PP2016} then proves the
claim for $C>0$. A fortiori it then also holds for all real $C$. $%
\blacksquare $

\bigskip

We are now ready to prove Theorem \ref{thm:pp24}. The proof follows the
structure of the proof of Theorem 5.1 in \cite{PP21}.

\textbf{Proof of Theorem \ref{thm:pp24}:} We apply Theorem \ref{thm:moregen}
with $m=n$ and $n_{j}=1$ for $j=1,\ldots ,m$, observing that then $\mathfrak{%
C}_{(n_{1},\ldots ,n_{m})}=\mathfrak{C}_{Het}$ and that condition (\ref%
{non-incl_Het_uncorr}) is equivalent to (\ref{non-incl_groupHet}) by Part
(b) of Lemma \ref{equivalence}. This then establishes that (\ref%
{size-control_Het}) follows from (\ref{non-incl_Het_uncorr}). The remaining
claim in Part 1 of Theorem \ref{thm:pp24} follows from Part 1 of Theorem \ref%
{thm:moregen}, if we can show that $\alpha ^{\ast }$ and $C^{\ast }$ given
in Theorem \ref{thm:moregen} can be written as claimed in Theorem \ref%
{thm:pp24}. To show this, we proceed as follows: Choose an element $\mu _{0}$
of $\mathfrak{M}_{0}$. Observe that $I_{1}(\mathcal{L}_{\#})\neq \emptyset $
(since $\dim (\mathcal{L}_{\#})<n-1<n$, cf.~Lemma \ref{lem:dimbound}), and
that for every $i\in I_{1}(\mathcal{L}_{\#})$ the linear space $\mathcal{S}%
_{i}=\func{span}(\Pi _{\mathcal{L}_{\#}^{\bot }}e_{i}(n))$ is $1$%
-dimensional (since $\mathcal{S}_{i}=\{0\}$ is impossible in view of $i\in
I_{1}(\mathcal{L}_{\#})$), and belongs to $\mathbb{J}(\mathcal{L}_{\#},%
\mathfrak{C}_{Het})$ in view of Proposition B.1 in Appendix B of \cite{PP21}
together with $\limfunc{dim}(\mathcal{L}_{\#})<n-1$. Since $T_{Het}$ is $G(%
\mathfrak{M}_{0})$-invariant (Remark C.1(i) in Appendix C of \cite{PP21}),
it follows that $T_{Het}$ is constant on $(\mu _{0}+\mathcal{S}%
_{i})\backslash \left\{ \mu _{0}\right\} $, cf.~the beginning of the proof
of Lemma 5.11 in \cite{PP3}. Hence, $\mathcal{S}_{i}$ belongs to $\mathbb{H}$
(defined in Lemma 5.11 in \cite{PP3}) and consequently for $C^{\ast }$ as
defined in that lemma%
\begin{equation}
C^{\ast }\geq \max \left\{ T_{Het}(\mu _{0}+\Pi _{\mathcal{L}_{\#}^{\bot
}}e_{i}(n)):i\in I_{1}(\mathcal{L}_{\#})\right\}  \label{ineq}
\end{equation}%
must hold (recall that $\Pi _{\mathcal{L}_{\#}^{\bot }}e_{i}(n)\neq 0$). To
prove the opposite inequality, let $\mathcal{S}$ be an arbitrary element of $%
\mathbb{H}$, i.e., $\mathcal{S}\in \mathbb{J}(\mathcal{L}_{\#},\mathfrak{C}%
_{Het})$ and $T_{Het}$ is $\lambda _{\mu _{0}+\mathcal{S}}$-almost
everywhere equal to a constant $C(\mathcal{S})$, say. Then Proposition B.1
in Appendix B of \cite{PP21} together with $\limfunc{dim}(\mathcal{L}%
_{\#})<n-1$ shows that $\mathcal{S}_{i}\subseteq \mathcal{S}$ holds for some 
$i\in I_{1}(\mathcal{L}_{\#})$. By Remark B.1(iv) given below, the condition
in (\ref{non-incl_Het_uncorr}) is equivalent to 
\begin{equation*}
e_{i}(n)\notin \mathsf{B}\text{ for every }i\in I_{1}(\mathcal{L}_{\#}).
\end{equation*}%
Therefore, (\ref{non-incl_Het_uncorr}) implies that we have $\mathcal{S}%
_{i}\nsubseteqq \mathsf{B}$ since $\Pi _{\mathcal{L}_{\#}^{\bot }}e_{i}(n)$
and $e_{i}(n)$ differ only by an element of $\mathcal{L}_{\#}$ and since $%
\mathsf{B}+\mathcal{L}_{\#}=\mathsf{B}$ (because of Part 1 of Lemma \ref%
{lem:additionV}). Thus $\mu _{0}+\mathcal{S}_{i}\nsubseteqq \mathsf{B}$ by
the same argument as $\mu _{0}\in \mathfrak{M}_{0}\subseteq \limfunc{span}%
(X) $ and $\mathsf{B}+\limfunc{span}(X)=\mathsf{B}$. We thus can find $s\in 
\mathcal{S}_{i}$ such that $\mu _{0}+s\notin \mathsf{B}$. Note that $s\neq 0$
must hold, since $\mu _{0}\in \mathfrak{M}_{0}\subseteq \limfunc{span}%
(X)\subseteq \mathsf{B}$. In particular, $T_{Het}$ is continuous at $\mu
_{0}+s$, since $\mu _{0}+s\notin \mathsf{B}$. Now, for every open ball $%
A_{\varepsilon }$ in $\mathbb{R}^{n}$ with center $s$ and radius $%
\varepsilon >0$ we can find an element $a_{\varepsilon }\in A_{\varepsilon
}\cap \mathcal{S}$ such that $T_{Het}(\mu _{0}+a_{\varepsilon })=C(\mathcal{S%
})$. Since $a_{\varepsilon }\rightarrow s$ for $\varepsilon \rightarrow 0$,
it follows that $C(\mathcal{S})=T_{Het}(\mu _{0}+s)$. Since $s\neq 0$ and
since $T_{Het}$ is constant on $(\mu _{0}+\mathcal{S}_{i})\backslash \left\{
\mu _{0}\right\} $ as shown before, we can conclude that $C(\mathcal{S}%
)=T_{Het}(\mu _{0}+s)=T_{Het}(\mu _{0}+\Pi _{\mathcal{L}_{\#}^{\bot
}}e_{i}(n))$, where we recall that $\Pi _{\mathcal{L}^{\bot }}e_{i}(n)\neq 0$%
. But this now, together with (\ref{ineq}), implies%
\begin{equation*}
C^{\ast }=\max \left\{ T_{Het}(\mu _{0}+\Pi _{\mathcal{L}_{\#}^{\bot
}}e_{i}(n)):i\in I_{1}(\mathcal{L}_{\#})\right\} .
\end{equation*}%
Using invariance of $T_{Het}$ w.r.t.~addition of elements of $\mathcal{L}%
_{\#}$ (cf.~Lemma \ref{lem:additionV}) we conclude that%
\begin{equation}
C^{\ast }=\max \left\{ T_{Het}(\mu _{0}+e_{i}(n)):i\in I_{1}(\mathcal{L}%
_{\#})\right\} .  \label{prelim}
\end{equation}%
Recall that $I_{1}(\mathcal{L}_{\#})\subseteq I_{1}(\mathfrak{M}_{0}^{lin})$%
. For $i\in I_{1}(\mathfrak{M}_{0}^{lin})\left\backslash I_{1}(\mathcal{L}%
_{\#})\right. $ we have $i\in I_{0}(\mathcal{L}_{\#})$, and thus $%
e_{i}(n)\in \mathcal{L}_{\#}$. Since $\mathcal{L}_{\#}\subseteq \mathsf{B}$, 
$e_{i}(n)\in \mathsf{B}$ follows. Using Part 1 of Lemma \ref{lem:additionV}
and $\mathfrak{M}_{0}\subseteq \mathsf{B}$, we conclude that $\mu
_{0}+e_{i}(n)\in \mathsf{B}$, and thus $T_{Het}(\mu _{0}+e_{i}(n))=0$. Since 
$T_{Het}$ is always nonnegative and since $I_{1}(\mathcal{L}_{\#})$ is
nonempty, we can write (\ref{prelim}) equivalently as%
\begin{equation*}
C^{\ast }=\max \left\{ T_{Het}(\mu _{0}+e_{i}(n)):i\in I_{1}(\mathfrak{M}%
_{0}^{lin})\right\} .
\end{equation*}%
The expression for $\alpha ^{\ast }$ given in the theorem now follows
immediately from the expression for $\alpha ^{\ast }$ given in Part 1 of
Theorem \ref{thm:moregen}.

Part 2-4 now follow from the corresponding parts of Theorem \ref{thm:moregen}
in light of what has been shown above. $\blacksquare $

\bigskip

\textbf{Remark B.1:} \emph{(Equivalent forms of the size-control conditions) 
}(i) The proof of Lemma \ref{equivalence} has shown that (\ref%
{non-incl_groupHet}) is not only equivalent to (\ref%
{non-incl_groupHet_uncorr}), but also to (\ref{non-incl_groupHetmod2a}) as
well as to (\ref{non-incl_groupHetmod2}).

(ii) Non-inclusion statements of the form "$\func{span}\left( \left\{
e_{i}(n):i\in J\right\} \right) \nsubseteqq \mathsf{B}$" ($J$ an index set)
appearing in (\ref{non-incl_groupHet}), (\ref{non-incl_groupHetmod2a}), and (%
\ref{non-incl_groupHetmod2}) can equivalently be written as "$e_{i}(n)\notin 
\mathsf{B}$ for some $i\in J$" due to the fact that $\mathsf{B}$ is a linear
space (as $R$ is $1\times k$). Similarly, "$\func{span}\left( \left\{
e_{i}(n):i\in J\right\} \right) \nsubseteqq \limfunc{span}(X)$" is
equivalent to "$e_{i}(n)\notin \limfunc{span}(X)$ for some $i\in J$".

(iii) In the special case where $m=n$ and $n_{1}=n_{2}=...=n_{m}=1$, we
learn from Lemma \ref{equivalence} and its proof that (\ref%
{non-incl_Het_uncorr}) is equivalent to (\ref{non-incl_groupHetmod2spez*}),
as well as to (\ref{delta}). Since $\mathcal{I}_{\#}^{c}\subseteq I_{1}(%
\mathcal{L}_{\#})\subseteq I_{1}(\mathfrak{M}_{0}^{lin})$ by Part 3 of Lemma %
\ref{gensol}, each one of (\ref{non-incl_Het_uncorr}), (\ref%
{non-incl_groupHetmod2spez*}), and (\ref{delta}) is in turn equivalent to
the condition%
\begin{equation}
e_{i}(n)\notin \limfunc{span}(X)\text{ for every }i\in I_{1}(\mathcal{L}%
_{\#}).  \label{simpl_1}
\end{equation}%
[As a point of interest we note that conditions (\ref{non-incl_Het_uncorr}),
(\ref{non-incl_groupHetmod2spez*}), (\ref{delta}), and (\ref{simpl_1}) are
in fact equivalent also if, in the notation of \cite{PP21}, we have $q\geq 1$%
, i.e., if a collection of $q$ restrictions is tested simultaneously. This
can be seen by an inspection of the proofs of these equivalences. However,
note that in case $q>1$ we have no result guaranteeing that \emph{these}
conditions are sufficient for size controllability of $T_{Het}$.]

(iv) Specializing Part (a) of Lemma \ref{equivalence} and its proof to the
case $n_{j}=1$ for $j=1,\ldots ,n=m$, and noting that $\mathcal{I}%
_{\#}^{c}\subseteq I_{1}(\mathcal{L}_{\#})$ (Lemma \ref{gensol}), one sees
that further equivalent forms of (\ref{non-incl_Het_uncorr}) are given by
the condition 
\begin{equation*}
e_{i}(n)\notin \mathsf{B}\text{ for every }i\in I_{1}(\mathcal{L}_{\#}),
\end{equation*}%
as well as by the condition%
\begin{equation*}
e_{i}(n)\notin \mathsf{B}\text{ for every }i\in \mathcal{I}_{\#}^{c},
\end{equation*}%
respectively. However, recall that while condition (\ref{non-incl_Het})
implies anyone of the two equivalent conditions above, it is, in general,
stronger in view of the examples in Appendix \ref{app:proofs}.

(v) Since in the special case where $m=n$ and $n_{1}=n_{2}=...=n_{m}=1$
condition (\ref{non-incl_Het_uncorr}) appears also as the size-control
condition for the standard (uncorrected) F-test statistic (see \cite{PP21}),
this condition can also be written in any of the equivalent forms given in
(iii) or (iv) in the case of testing a single restriction as considered
here. [The equivalence of (\ref{non-incl_Het_uncorr}) with the other
conditions in (iii) above even holds in the more general case where more
than one restriction is subject to test.] We note that the before given
equivalences do \emph{not} rely on Assumption \ref{R_and_X}, an assumption
that also does not appear in the size control results in \cite{PP21} for the
classical (uncorrected) F-test statistic.

\textbf{Remark B.2:} The proof of Theorem \ref{thm:pp24} shows that $C^{\ast
}$ defined in (\ref{Cstar}) can alternatively be written as in (\ref{prelim}%
). The representation (\ref{prelim}) has two advantages over (\ref{Cstar}):
First, the index set $I_{1}(\mathcal{L}_{\#})$ is potentially smaller than $%
I_{1}(\mathfrak{M}_{0}^{lin})$ (see Lemma \ref{gensol}); second, since $%
e_{i}(n)\notin \mathsf{B}$ for $i\in I_{1}(\mathcal{L}_{\#})$ under
condition (\ref{non-incl_Het_uncorr}) (see Remark B.1(iv)), also $\mu
_{0}+e_{i}(n)\notin \mathsf{B}$ for such $i$ ($\mu _{0}\in \mathfrak{M}_{0}$%
). Thus, (\ref{prelim}) does not rely on the way $T_{Het}$ has been defined
on the set $\mathsf{B}$.

\bibliographystyle{ims}
\bibliography{refs}

\end{document}